\documentclass[11pt]{article}

\usepackage{amsmath}
\usepackage{amssymb}
\usepackage{amscd}
\begin{document}
\newtheorem{theorem}{Theorem}[section]
\newtheorem{lemma}[theorem]{Lemma}
\newtheorem{corollary}[theorem]{Corollary}
\newtheorem{definition}[theorem]{Definition}
\newtheorem{proposition}[theorem]{Proposition}
\newtheorem{theo}[theorem]{Theorem}
\newtheorem{prop}[theorem]{Proposition}
\newtheorem{coro}[theorem]{Corollary}
\newtheorem{lemm}[theorem]{Lemma}
\newtheorem{defi}[theorem]{Definition}
\newtheorem{rem}[theorem]{Remark}
\newtheorem{remark}[theorem]{Remark}
\newtheorem{ex}[theorem]{Example}
\newtheorem{example}[theorem]{Example}
\renewcommand{\theequation}{\thesection.\arabic{equation}}
\def\arp{\leftharpoonup}
\def\iso{\widetilde {\longrightarrow}}
\def\B{\mathcal B}
\def\qsl{\mathcal O(SL_q(2))}
\def\fre{k[z,z^{-1}] * \mathcal O(SL_q(2))}
\def\N{\mathbb N}
\def\Z{\mathbb Z}
\def\CoInn{{\rm CoInn}}
\def\CoInt{{\rm CoInt}}
\def\CoOut{{\rm CoOut}}
\def\C{\mathbb C}
\def\R{\mathbb R}
\def\Reg{{\rm Reg}}
\def\square{\Box}
\def\argh{\rightharpoonup}
\def\m#1{a_{#1}}
\def\g#1{g_{#1}}
\def\h#1{h_{#1}}
\def\k#1{b_{#1}}
\def\M{{\cal M}}
\def\Z{{\mathbb Z}}
\def\F{{\cal F}}
\def\mm#1{m_{#1*}}
\def\a#1{a_{#1}}
\def\aa#1{a_{#1*}}
\def\b#1{b_{#1}}
\def\bb#1{b_{#1*}}
\def\l#1{l_{#1}}
\def\ad{{\rm ad}}
\def\n#1{n_{#1}}
\def\nn#1{n_{#1*}}
\def\totimes{\tilde{\otimes}}
\def\c#1{c_{#1}}
\def\cc#1{c_{#1*}}
\def\coh#1{H^{#1}_L(A)}
\def\id{{\rm id}}
\def\ydh{{\cal YD}(A)}
\def\yddh{{\cal YD}(A^*)}
\def\lefth{_{A}{\cal M}}
\def\lefthd{_{A^*}{\cal M}}
\def\righth{{\cal M}^A}
\def\righthd{{\cal M}^{A^*}}
\def\CC{{\cal C}}
\def\D{{\cal D}}
\def\RR{{\cal R}}
\def\U{{\cal U}}
\def\rsb{r_{s,\beta}}
\def\rsg{r_{s,\gamma}}
\def\f{{\bf k}}
\def\km{\f^{*}}
\def\lazy{Z^2_L(A)}
\def\pf{\noindent{\bf Proof:} }
\def\cll{\Reg_{aL}^1(A)}
\def\reg#1{\Reg^#1(A)}
\def\regl#1{\Reg_L^#1(A)}
\def\bl{B^2_{aL}(A)}
\def\aut{{\rm Aut}_{\rm Hopf}(A)}
\def\lazyn{Z_L^n(A)}
\def\lu{\rightharpoonup}
\def\bic{B \! \Join \! A}
\def\cp{\mathcal{ZP}}
\def\lu{\rightharpoonup}
\def\ru{\leftharpoonup}

\title{\textbf{\textsc{Lazy cohomology:
an analogue of the Schur multiplier for arbitrary
Hopf algebras}}} 
\author{\textsc{Julien Bichon$^*$ and Giovanna Carnovale$^{**}$}}
\date{\small{\textsl{*Laboratoire de Math\'ematiques Appliqu\'ees,\\
Universit\'e de Pau et des Pays de l'Adour, \\
IPRA, Avenue de l'universit\'e,
64000 Pau, France.}\\
E-mail: Julien.Bichon@univ-pau.fr}\\
~\\
{\textsl{{**}Dipartimento di Matematica Pura ed Applicata,\\ 
Universit\`a di Padova,\\
via Belzoni 7, 35131 Padova, Italy.}\\
E-mail: carnoval@math.unipd.it}\\}

\makeatletter
\renewcommand{\@makefnmark}{}
\makeatother

%\begin{document}

\maketitle

\begin{abstract}
We propose a detailed systematic study of a group 
$H^2_L(A)$ associated, by elementary means of lazy 2-cocycles,
to any Hopf algebra $A$. This group was introduced
by Schauenburg in order to generalize G.I. Kac's exact sequence.
We study the various interplays of lazy cohomology
in Hopf algebra theory: Galois and biGalois objects, 
Brauer groups and projective representations.
We obtain a Kac-Schauenburg-type sequence
for double crossed products of possibly infinite-dimensional 
Hopf algebras. Finally the explicit computation of $H^2_L(A)$ 
for monomial Hopf algebras and for a class of
cotriangular Hopf algebras is performed.
%a group $H_L^2(A)$ to any Hopf algebra
%$A$. 
%This group may reasonably be seen as an analogue
%of the Schur multiplier, since when $A = k[G]$
%is a group algebra, we have $H_L^2(k[G]) = H^2(G, k^\cdot)$.
%The Hopf-Galois extensions corresponding to
%lazy 2-cocycles are characterized
%as Hopf-Galois extensions with a special symmetry property. The group
%$H_L^2(A)$ is embedded as a normal subgroup
%of BiGal$(A)$, the biGalois group of $A$ and a group
%exact sequence relating the semi-direct product
%${\rm CoOut}(A) \ltimes H_L^2(A)$ and ${\rm BiGal}(A)$ is constructed.
%A Schur-Yamazaki type formula for tensor products is also provided.
%Finally, the group $H_L^2(A)$ enables us to define a 
%monoidal $H^2_L(A)$-category of projective representations of a Hopf 
%algebra.
%with a structure closely related to the notion of $\pi$-category
%introduced by Turaev in homotopy quantum field theory.
%The explicit computation of $H^2_L(A)$ 
%for monomial Hopf algebras and for a class of
%cotriangular Hopf algebras is performed.
%Finally, using the classification of finite dimensional triangular
%$Hopf algebras over ${\mathbb C}$ we are able to compute $H^2$ for 
%finite-dimensional cotriangular Hopf algebras over $\mathbb C$. 
\end{abstract}



\section*{Introduction}

In 1968 Sweedler defined in \cite{[Sw]} the cohomology for a 
cocommutative Hopf
algebra $H$ with coefficients in a module algebra $A$,  
relating it to the Brauer group of a field and to well known
cohomology theories such as Lie algebra cohomology, group cohomology, Galois
cohomology. For instance, he showed that when the algebra $A$ is
commutative, there is a bijection between the second cohomology group
and the equivalence classes of $H$-cleft extensions of $A$. 

After this important achievement, several results have
been obtained in the direction of providing a
cohomological interpretation of many
natural constructions arising in Hopf algebra theory in the 
non-cocommutative case. For instance, the theory of cleft Hopf-Galois 
extensions can be described in terms of $2$-cocycles up
to coboundaries (\cite{doitake1}, \cite{doi}, \cite{bcm}). However, the 
words
``cocycle'' and ``coboundary'' are used because they
satisfy equations that look like Sweedler's conditions and coincide
with it in the cocommutative case but they do not really correspond 
to a given complex. 

The basic problem with Hopf 2-cocycles is that 
the convolution product of two 2-cocycles
is no longer a  2-cocycle in general.
This problem disappears when we deal with 
what we call lazy 2-cocycles, i.e. those cocycles that
are convolution commuting with the product of $A$. In this
way one can associate a group $H^2_L(A)$ to any Hopf algebra 
$A$ with $H^2_L(k[G]) =H^2(G,k^\cdot)$ in the case of a group algebra,
so it might reasonably seen as an analogue of the Schur multiplier
for Hopf algebras. 
The group $H^2_L(A)$ was introduced by Schauenburg in \cite{[Sc3]}
(under the notation $H^2_c(A)$) in his generalization of G. I.
Kac's exact sequence \cite{[Kac]}.
The Kac exact sequence is a useful tool for
computing Hopf algebras extension by group algebras.

\medskip

In this paper we propose a detailed systematic study
of the group $H^2_L(A)$, that we call the (second)
lazy cohomology group of $A$. Apart from the 
Kac-Schauenburg exact sequence,  of which we give a 
generalization in the case of double crossed product
Hopf algebras, to possibly infinite-dimensional Hopf
algebras, 
the motivation for this study has several origins, one of which is 
the study of the biGalois group of $A$.
This group, denoted by BiGal$(A)$,
 might be defined as the group of isomorphism
classes of linear monoidal auto-equivalences of the category
of $A$-comodules. According to Schauenburg \cite{[Sc1]}, this group
might also be described as the set of isomorphism classes
of $A$-$A$-biGalois extensions endowed with the cotensor product
(we assume for simplicity in this introduction
that the base ring is a field). Schauenburg's description
is certainly the most efficient one for concrete computations
of BiGal$(A)$. The simplest example is when $A= k[G]$ is a group
algebra and we have BiGal$(k[G]) = {\rm Aut}(G) \ltimes 
H^2(G,k^\cdot)$.
However it is in general a difficult task to give an explicit
description of the biGalois group, even when Gal$(A)$, 
the set of isomorphism classes of right $A$-Galois objects, 
has been determined. For example Schauenburg \cite{[Sc2]}
has described the set Gal$(A\otimes B)$ for a Hopf algebra
tensor product using Gal$(A)$ and Gal$(B)$ but a similar 
description at the biGalois group level seems to be still unknown.
Schauenburg's computation of the biGalois groups
of the generalized Taft algebras $H_{N,m}$ with $m$ grouplike
also shows that the complete description of the
biGalois group is a more delicate problem, in general, than
the description of the Galois objects. Some simplifications
of Schauenburg's formula for $H_{N,m}$,
at least from the theoretical viewpoint, are given
in \cite{[BiG]} but there is certainly still more
to be understood.

The group $H_L^2(A)$ is studied in this perspective.
It may be realized as a normal subgroup of BiGal$(A)$, and
hence might serve as a first approximation to understand
the structure of  BiGal$(A)$.
From the monoidal categories viewpoint, it is 
the subgroup of BiGal$(A)$ consisting of isomorphism
classes of linear monoidal auto-equivalences of the category
of $A$-comodules that are isomorphic, \textsl{as functors},
to the identity functor.
The technical simplification with $H_L^2(A)$ is that
once Gal$(A)$ has been determined, it is much easier 
to describe than  BiGal$(A)$. 
A group morphism from \CoOut$(A)\ltimes H_L^2(A)$ to BiGal$(A)$
is constructed (this morphism is an isomorphism when $A$ is
cocommutative \cite{[Sc1]}), and the corresponding exact sequence 
is described.
Coming back to Hopf tensor products, a Schur-Yamazaki type formula,
derived from the generalized Kac-Schauenburg exact sequence, 
is given, describing $H_L^2(A \otimes B)$ from $H_L^2(A)$,
$H_L^2(B)$ and the group of central $A\otimes B$ pairings.
This formula generalizes the classical one (see e.g. \cite{[Kar]})
describing the second cohomology group of a direct product of groups.
 
\smallskip

Another motivation for the study of lazy cohomology comes from the
theory of Brauer groups of a Hopf algebra.
Although a systematic description was not given yet, 
lazy cocycles were used extensively for the computation of the Brauer
groups $BC(k, H, r)$ (see \cite{CVOZ} for the construction of $BC$)  
of the families of coquasitriangular Hopf algebras $E(n)$
and $H_\nu$ (see \cite{gio} where the results in \cite{yinhuo} are 
generalized,
\cite{GioJuan1}, \cite{GioJuan3}). According to the known examples, we
 expect that the lazy cohomology
group will occur as a subgroup of the Brauer group of any cotriangular
Hopf algebra.

\smallskip

Finally, although there exists a nice cohomology theory for Hopf
algebras, with coefficients into Hopf bimodules (see \cite{taillefer}
and the references therein), there is not in this
framework a cohomology with coefficients in the multiplicative group
of the field, which is used for example in projective representation
theory
of groups. Our group $H^2_L(A)$ enables us to construct a monoidal
$H^2_L(A)$-category of projective representations of $A$. Such 
categorical structures were
considered by Turaev in homotopy quantum field theory (\cite{[Tu]}).

\medskip

This paper is organized as follows.
The first Section is devoted to the direct elementary 
construction of $H_L^2(A)$ using lazy 2-cocycles, and 
to other preliminary constructions. Section \ref{exa} contains
the first few direct computations with explicit cocycles for Sweedler's
Hopf algebra $H_4$ and for $E(n)$. 
The setup of Section \ref{galois} is the very classical link
between Hopf-Galois extensions and $2$-cocycles: 
the Galois objects corresponding to lazy cocycles
are characterized by a special symmetry property. The lazy cohomology
group is then embedded as a normal subgroup of the group of BiGalois
objects and it is shown to be cocycle twist-invariant. 
In Section \ref{KS} we formulate and prove a  Kac-Schauenburg-type sequence
for double crossed product Hopf algebras, generalizing it
to the possibly infinite-dimensional case.
We derive from it a Schur-Yamazaki-type formula so that
the description of the lazy cohomology for the Drinfeld double of a
(finite-dimensional) Hopf algebra is
given as a corollary. Section \ref{universal} deals with the connection between lazy
cohomology and universal $R$-forms for coquasitriangular Hopf
algebras with applications to the Brauer group of a Hopf algebra. 
In Section \ref{projective} we introduce a monoidal category of projective
representations for Hopf algebras. This category has the structure of monoidal
$H^2_L(A)$-category over the base field $k$. 
The last two sections deal with examples: 
in Section \ref{monomial} the lazy
cohomology group of a monomial Hopf algebra is computed. Section
\ref{special} 
contains the description of lazy cohomology for those
finite-dimensional cotriangular Hopf algebras for which the linear action of
the group datum in Andruskiewitsch, Etingof and Gelaki's terminology  
is faithful. An Appendix deals with relations between
lazy cocycles and general 
crossed systems.

\subsection*{Notation and conventions}

Unless otherwise stated $k$ is a commutative ring. Unadorned tensor
products will be over $k$. The group of invertible elements in
$k$ will be denoted by $k^\cdot$. For any pair of $k$-modules $V$ and
$W$, the usual  flip map $V\otimes W\to W\otimes V$ interchanging the
tensorands will be denoted by $\tau$. For a Hopf algebra over $k$
the  product will be denoted by $m$,
the coproduct by $\Delta$ and the antipode by $S$. We adopt a Sweedler's
like notation e.g.: $\Delta(a)=\a1\otimes\a2$.
  
\section{Lazy cocycles} 

Let $A$ be a Hopf algebra. In this section we give the 
detailed elementary 
direct construction, using lazy cocycles,
of the group $H_L^2(A)$ and some further preliminary constructions.

\smallskip 

We first recall some classical notions.
The set of convolution invertible linear maps 
$\mu : A \to k$ satisfying $\mu(1)=1$ is denoted
by $\Reg^1(A)$. Similarly  the set of convolution invertible linear 
maps 
$\sigma : A \otimes A \longrightarrow k$ satisfying 
$\sigma(a,1)=\varepsilon(a) = \sigma(1,a)$, for all $a \in A$, is 
denoted
by $\Reg^2(A)$. It is clear that both $\Reg^1(A)$ and $\Reg^2(A)$ are 
groups
under the convolution product.

\begin{definition}
An element $\mu \in \Reg^1(A)$ is said to be \textbf{lazy}
if $\mu * {\rm id}_A= {\rm id}_A * \mu$. The set of lazy
elements of $\Reg^1(A)$, denoted $\Reg^1_L(A)$, is a central subgroup 
of 
$\Reg^1(A)$.

An element $\sigma \in \Reg^2(A)$ is said to be \textbf{lazy}
if $\sigma * m= m * \sigma$, which is to say that
\begin{equation}\label{condition}\sigma(a_1,b_1)a_2b_2 = 
\sigma(a_2,b_2)a_1b_1, \ \forall a,b \in A.\end{equation}
The set of lazy
elements of $\Reg^2(A)$, denoted $\Reg_L^2(A)$, is a subgroup of 
$\Reg^2(A)$.
\end{definition}

 If $A$ is cocommutative then
$\Reg^q_L(A)=\Reg^q(A)$ for $q=1,2$.

We recall that a \textbf{ left 2-cocycle} over $A$ is an element
$\sigma \in \Reg^2(A)$ such that 
$$\sigma(a_{1}, b_{1}) \sigma(a_{2}b_{2},c) =
\sigma(b_{1},c_{1}) \sigma(a,b_{2} c_{2}), \ \forall a,b,c \in A.$$

We shall denote by $Z^2(A)$ the set of left $2$-cocycles and by
$Z^2_L(A)$ the set $Z^2(A)\cap \Reg^2_L(A)$ of {\bf lazy $2$-cocycles}.
 
Similarly a \textbf{right 2-cocycle} over $A$ is an element
$\sigma \in \Reg^2(A)$ such that 
$$\sigma(a_{1}b_{1},c) \sigma(a_{2},b_{2}) =
\sigma(a,b_{1}c_{1}) \sigma(b_{2}, c_{2}), 
\ \forall a,b,c \in A.$$
It is well known that
$\sigma\in {\rm Reg}^2(A)$ is a left 2-cocycle if and only if 
$\sigma^{-1}$ is a right 2-cocycle.

The basic problem with left (or right) 2-cocycles is that 
the convolution product of two left 2-cocycles
is no longer a left 2-cocycle in general.
This problem disappears when we deal with lazy 2-cocycles.  

\begin{lemm}\label{group} Let $\sigma,\omega\in Z^2(A)$.
\begin{enumerate}
\item If $\omega \in Z_L^2(A)$, then $\sigma *\omega\in Z^2(A)$.
\item If $\sigma \in Z_L^2(A)$ then $\sigma$ is a right 2-cocycle
and $\sigma^{-1}$ is a left 2-cocycle.
\end{enumerate} 
In particular, $Z_L^2(A)$ is a subgroup of $\Reg_L^2(A)$.
\end{lemm}
\pf The proof is by straightforward computation. We shall perform it
in order to show how condition (\ref{condition}) can be used.

\noindent{\it 1.} Let $\omega\in Z_L^2(A)$, let $\sigma\in Z^2(A)$ and
let $a,b,c\in A$. Then
$$
\begin{array}{rl}
\sigma*\omega(\a1,\b1)\sigma*\omega(\a2\b2,c)&=\sigma(\a1,\b1)\sigma(\a2\b2,\c1)\omega(\a3,\b3)\omega(\a4\b4,\c2)\\
&=\sigma(\b1,\c1)\sigma(\a1,\b2\c2)\omega(\b3,\c3)\omega(\a2,\b4\c4)\\
&=\sigma(\b1,\c1)\omega(\b2,\c2)\sigma(\a1,\b3\c3)\omega(\a2,\b4\c4)\\
&=\sigma*\omega(\b1,\c1)\sigma*\omega(\a1,\b2\c2).
\end{array}
$$ 

\noindent{\it 2.} Let $\sigma\in Z^2_L(A)$ and let $a,\,b,\,c\in H$. 
Then 
$$
\begin{array}{rl}
\sigma(\b2,\c2)\sigma(a, \b1\c1)&=
\sigma(\b1,\c1)\sigma(a, \b2\c2)\\
&=\sigma(\a1,\b1)\sigma(\a2\b2, c)\\
&=\sigma(\a1\b1, c)\sigma(\a2,\b2)
\end{array}
$$
where for the first and the third equality we used that $\sigma$ is
lazy and for the second equality we used that $\sigma$ is a left
2-cocycle.\hfill$\Box$  
%%%%%It is well known that $\sigma$ is a right $2$-cocycle if
%and only if $\sigma^{-1}$ is a left $2$-cocycle. 

\bigskip

The terminology {\em lazy} for $2$-cocycles
is motivated by the fact that a lazy $2$-cocycle
does not alter the  Hopf algebra $A$ through Doi's twisting 
procedure, i.e., if $\sigma$ is a lazy $2$-cocycle,
 $\id_A$ is a Hopf algebra isomorphism between $A$ and the Hopf
algebra $_\sigma A_{\sigma^{-1}}$ that has the underlying
coalgebra as $A$ and
product given by 
\begin{equation}a\cdot_\sigma 
b=\sigma(\a1,\b1)\a2\b2\sigma^{-1}(\a3,\b3).
\end{equation}

Similarly, given $\gamma\in \Reg^1(A)$, 
the measuring $\ad(\gamma)\colon A\to {\rm End}(A)$ given by
$\ad(\gamma)(a)=\gamma^{-1}(\a1)\a2\gamma(\a3)$ is trivial if and only 
if
$\gamma$ is lazy. 

\begin{example}{\rm It is well known that if $A$ is coquasitriangular
with
universal $r$-form $r$, then $r$ is a left $2$-cocycle for
$A$. The form $r$ is lazy if and only if $A$ is commutative.}
\end{example}

\begin{example}{\rm Let $A$ be a coquasitriangular Hopf algebra and let 
$r,\,s$ lie in
the set of coquasitriangular structures $\U$ of $A$. Then
$(r\circ\tau)*s\in \lazy$. Indeed for $a,b,c\in A$ we have: 
%Then, since $s\in\U$ we have
$$
\begin{array}{l}
r(\m1,\k1)s(\k2,\m2)r(\k3\m3,\c1)s(\c2, \k4\m4)\\
=r(\m1,\k1)r(\m2\k2,\c1)s(\k3,\m3)s(\c2, \k4\m4).
\end{array}
$$
Since universal $r$-forms are $2$-cocycles the above is equal to:
$$
\begin{array}{l}
r(\k1,\c1)r(\m1,\k2\c2)s(\c3,\k3)s(\c4\k4, \m2)\\
=r(\k1,\c1)r(\m1,\c3\k3)s(\c2,\k2)s(\c4\k4, \m2)\\
=(r\tau *s)(\c1,\k1)(r\tau*s)(\c2\k2, a)
\end{array}
$$
where we used again that $s\in\U$. Hence $(r\circ\tau)* s$ is a 
$2$-cocycle.
It is not hard to see that $(r\circ\tau)*s$ satisfies condition
(\ref{condition}). It follows that $Z_L^2(A)$ is nontrivial for 
coquasitriangular Hopf algebras which are not cotriangular.}
\end{example}

\begin{remark}\label{crossed}{\rm By part 1 of Lemma \ref{group}, if
  $\sigma$ is any $2$-cocycle and
$\omega\in\lazy$ then $\sigma*\omega$ is a $2$-cocycle. Similarly if
$\sigma$ is a right $2$-cocycle and $\omega\in\lazy$ then
$\omega*\sigma$ is a right $2$-cocycle. 
Hence $\lazy$ acts on the right on $Z^2(A)$ by $\sigma\mapsto
\sigma*\omega$. More generally, one can define a right action of
$Z^2_L(A)$ on the set of general crossed systems. This will be
done in the Appendix.}
\end{remark}

The group $H_L^2(A)$ will be defined as a quotient group
of $Z_L^2(A)$. For this we need a ``differential'' connecting
$\Reg_L^1(A)$ and $\Reg_L^2(A)$. In order to do so, we simply have to 
mimic what happens in the cocommutative case (see \cite{[Sw]}).

The map
$$\partial : \Reg^1(A) \longrightarrow \Reg^2(A)$$
is defined by $\partial(\mu) = (\mu\otimes \mu) * (\mu^{-1} \circ m)$
for $\mu \in \Reg^2(A)$, that is: $\partial(\mu)(a,b) =
\mu(a_1)\mu(b_1) \mu^{-1}(a_2b_2)$ for $a,b \in A$.
Some basic properties of the operator $\partial$ are given 
in the following lemma.

\begin{lemm}\label{homom} Let $\mu, \phi \in \Reg^1(A)$.
\begin{enumerate}
\item $\partial(\mu) = \varepsilon \otimes \varepsilon$
if and only if $\mu \in {\rm Alg}(A,k)$.
\item If $\mu \in \Reg_L^1(A)$, then $\partial(\mu) \in \Reg_L^2(A)$.
\item $\partial(\mu*\phi) = (\mu \otimes \mu) * \partial(\phi) * 
(\mu^{-1}\circ m)$.
\item If $\mu \in \Reg_L^1(A)$, then $\partial(\mu*\phi) = 
\partial(\phi) *
\partial (\mu)$.
\item If $\mu \in \Reg_L^1(A)$, then for 
$\sigma \in \Reg_L^2(A)$, we have $\partial(\mu)*\sigma = \sigma *
\partial(\mu)$. 
\item If $\partial(\phi) \in \Reg_L^2(A)$, then 
$\partial(\mu*\phi) = \partial(\mu) *
\partial (\phi)$.
\item$\partial(\mu)$ is a left 2-cocycle.
\end{enumerate}
In particular the map $\partial$ induces 
a group morphism $\Reg_L^1(A) \to Z_L^2(A)$ with image contained in the
center of $Z_L^2(A)$. 
\end{lemm}
\pf The proof of all statements follows by direct computation and we
leave it to the reader.\hfill$\Box$

\medskip

Lemma \ref{homom} leads to the following definition:

\begin{defi}
Let $A$ be a Hopf algebra. The lazy cohomology groups
$H_L^1(A)$ and $H_L^2(A)$ are defined in the following way:
$$H_L^1(A) := {\rm Ker}(\partial_{|\Reg_L^1(A)})
= \{\mu \in {\rm Alg}(A,k)~|~\mu * {\rm id}_A = {\rm id}_A * \mu \},$$
$$H_L^2(A) := Z_L^2(A)/ B_L^2(A),$$
where $B_L^2(A)$ is the central subgroup 
$\partial(\Reg_L^1(A))$ of $Z_L^2(A)$.
\end{defi}

Let us observe that, according to Lemma \ref{homom}, 
lazy cocycles belonging to the same class in
$H^2_L(A)$ are cohomologous in the the sense of \cite{doi}.

The group $H_L^1(A)$ is obviously commutative. Although all the
examples of $H^2_L(A)$ computed here are commutative,
we see no reason why this group should be commutative in general. 

\begin{rem}\label{coquasi}
{\rm Schauenburg has defined more generally the lazy cohomology group,
denoted by ${\cal H}^2_c(H)$ where $c$ stands for central,  
of a coquasibialgebra $H$ in \cite[Section 6]{[Sc3]}. This construction is in fact very natural
in view of the monoidal category interpretation 
of lazy cohomology as given further in Remark \ref{auto-equivalence}.}
\end{rem}

\begin{rem}\label{absolute}
{\rm We could also define an abelian analogue 
of the Schur multiplier in the following way. 
Let us say that a 2-cocycle $\sigma$ is absolutely
central if 
$$\sigma(a_1,b_1) a_2 \otimes b_2 = 
\sigma(a_2,b_2) a_1 \otimes b_1, \ \forall a,b \in A.$$ 
The absolutely central 2-cocycles clearly form a central 
subgroup of $\Reg_L^2(A)$, denoted by $Z_{\rm ab}^2(A)$.
We do not necessarily have $\partial(\Reg_L^1(A))\subset Z_{\rm 
ab}^2(A)$ but
we still can define an abelian group $H^2_{\rm ab}(A)$ as the
subgroup  of $H_L^2(A)$ generated by the classes of absolutely central
2-cocyles. We focus on the present $H_L^2(A)$ for several
reasons: absolute centrality
seems to be very
restrictive and already for Sweedler's Hopf algebra the set of
absolutely central $2$-cocycles is trivial; 
we wanted our analogue of the Schur multiplier
to be the biggest possible as a subgroup of ${\rm BiGal}(A)$;
there is a nice property corresponding to 
the notion of lazy cocycle at the Hopf-Galois
extension level (see Section \ref{galois}) while we have not found
such a property  for absolutely central 2-cocycles; lazy cocyles
turned out to be a useful tool in the computation of the Brauer group 
of a 
coquasitriangular  Hopf algebra. From the known examples (see
\cite{gio}, \cite{GioJuan1} and \cite{GioJuan3}) 
we expect that lazy cohomology would occur as a subgroup of the Brauer
group of any cotriangular Hopf algebra.}
\end{rem}

\begin{rem}{\rm It is possible to generalize the notion of laziness
to higher degree cochains and to adopt a suitable form of 
Sweedler's operators $D^q$. In particular, one may view
the lazy $2$-cocycles as those cochains $\omega$ for which a suitable 
variation
of $D^2$ gives $D^2\omega=\varepsilon^{\otimes3}$. 
However, even if for degree $3$ we
still have the inclusion $B^3_L(A)\subset Z^3_L(A)$, at a first sight
there seems to be no natural group structure on $Z^3_L(A)$ nor on 
the quotient $H_L^3(A)$.}
\end{rem}

\medskip

We have seen that $\partial(\Reg^1_L(A))\subset Z^2_L(A)$. However, it
might happen that if $\gamma\in \Reg^1(A)$ and $\gamma\not\in
\Reg^1_L(A)$ still $\partial(\gamma)\in Z^2_L(A)$. We shall analyze
such elements $\gamma$. 

Let $\gamma\in \Reg^1(A)$, let $\sigma\in Z^2(A)$, and let
$\sigma^\gamma$ be the $2$-cocycle, cohomologous to $\sigma$ in the 
usual sense defined by:
$\sigma^\gamma
:=(\gamma\otimes\gamma)*\sigma*(\gamma^{-1}\circ m)$.
It is well known that cohomologous cocycles yield isomorphic Hopf
algebra twists $_\sigma\!A_{\sigma^{-1}}$ and
$_{\sigma^\gamma}\!A_{(\sigma^{\gamma})^{-1}}$ with isomorphism
$\ad(\gamma)\colon\; _\sigma\!A_{\sigma^{-1}}\to
 _{\sigma^\gamma}\!A_{(\sigma^{\gamma})^{-1}}$ given by
$\gamma^{-1}*\id_A*\gamma$. In particular, if 
$\sigma=\varepsilon\otimes\varepsilon$ then
$\ad(\gamma)$ 
is an Hopf algebra isomorphism
$A\to\,_{\partial(\gamma)}\!A_{(\partial(\gamma))^{-1}}$. 
Therefore, if
$\partial(\gamma)\in Z^2_L(A)$, $\ad(\gamma)\in\aut$, the group of
Hopf algebra
automorphisms of $A$. In fact, condition $(\ref{condition})$ for
$\partial(\gamma)$ is equivalent to the requirement that 
$\ad(\gamma)$ is an Hopf algebra
automorphism. We have proved the following
\begin{lemma}A coboundary $\partial(\gamma)\in \Reg^2_L(A)$  if and 
only
if $\ad(\gamma)$ is an Hopf algebra automorphism.  \hfill$\Box$
\end{lemma}

We call an Hopf automorphism of type $\ad(\gamma)$ a {\bf cointernal}
automorphism. We shall denote the set of cointernal automorphisms of 
$A$
by $\CoInt(A)$. 

\begin{lemma}\label{list}The following properties hold:
\begin{enumerate}
\item $\ad\colon \Reg^1(A)\to{\rm Aut}_{\rm coalg}(A)$ is a group
  morphism with kernel $\Reg^1_L(A)$.
\item The set
$$
\begin{array}{rl}
\cll:&=\{\gamma\in \Reg^1(A)~|~\partial(\gamma)\in \Reg^2_L(A)\}\\
&=\{\gamma\in \Reg^1(A)~|~\ad(\gamma)\in \CoInt(A)\}\\
&=\ad^{-1}(\aut)
\end{array}
$$
is a subgroup of $\reg1$. An element of $\cll$ is called almost
lazy. 
\item ${\rm CoInt}(A)$ is a normal subgroup of $\aut$.
\item The normal subgroup  $\CoInn(A)$ of $\aut$ given by:
$$
\CoInn(A)=
\{f\in \aut\;|\;\exists \phi\in {\rm Alg}(A,\,k) \mbox{ with }
f=(\phi\circ S)*\id_A* \phi\}
$$
is contained in $\CoInt(A)$. 
\item $\partial\colon \cll\to Z^2_L(A)$ is a group morphism and its 
kernel
is  ${\rm Alg}(A,\,k)$. 
\end{enumerate} 
\end{lemma}
\pf It is well known that
$\ad(\gamma*\theta)=\ad(\gamma)\circ\ad(\theta)$ for every $\gamma$
and $\theta\in \Reg^1(A)$. This implies most of the
statements. Statement $5$ follows from part 6 of Lemma
\ref{homom}.
\phantom{.}\hfill$\Box$ 

\medskip 

The knowledge of $\partial(\cll)=B^2(A)\cap Z^2_L(A)$
helps us 
 to detect when $\CoInn(A)=\CoInt(A)$.

\begin{lemma}\label{equal}With notation as before,
$B^2(A)\cap \Reg^2_L(A)=B^2_L(A)$ if and only if $\CoInn(A)=\CoInt(A)$.
\end{lemma}
\pf Let us suppose that $B^2(A)\cap \Reg^2_L(A)=B^2_L(A)$. For every 
$\gamma\in\cll$ there
exists $\chi\in\regl1={\rm Ker}(\ad)$ such that 
$\partial(\gamma)=\partial(\chi)$. Then
$\partial(\gamma*\chi^{-1})=1$, i.e., $\gamma*\chi^{-1}$ is an 
algebra morphism $A\to k$. Then, 
$$\ad(\gamma)=\ad(\gamma)\circ\ad(\chi)^{-1}=\ad(\gamma*\chi^{-1})\in
\CoInn(A).$$  

Let us now suppose that $\CoInn(A)=\CoInt(A)$ and let $\sigma$ be any
element of $\partial(\cll)$. Then $\sigma=\partial(\gamma)$ for some
$\gamma\in\cll$ and $\ad(\gamma)=\ad(\chi)$ for some 
algebra map $A\to k$ so that
$\partial(\chi^{-1})=\partial(\chi\circ 
S)=\varepsilon\otimes\varepsilon$ and
$\chi\circ S\in\cll$. 

Therefore, $\ad({\gamma*(\chi\circ S)})=\id_A$ so that
$\gamma*(\chi\circ S)\in\regl1$ and 
$$\sigma=\partial(\gamma)=\partial(\gamma)*\partial(\chi\circ 
S)=\partial(\gamma*(\chi\circ S))\in B^2_L(A).$$ 
\hfill$\Box$        

\medskip

We define the following group:
\begin{defi}
Let $A$ be a Hopf algebra. 
The group ${\rm CoOut}^-(A)$ is
$${\rm CoOut}^-(A) :=
\cll/
{\rm ad}^{-1}({\rm \CoInn}(A))\cong \CoInt(A)/\CoInn(A).$$
\end{defi}

By the second equality ${\rm CoOut}^-(A)$ can be viewed as a subgroup 
of ${\rm
  CoOut(A)}=\aut/\CoInn(A)$, the
group of co-outer automorphisms of $A$. In Example \ref{counterex} we 
shall consider a Hopf algebra for which ${\rm CoOut}^-(A)$ is not
  trivial. 

\medskip

%%%%%This terminology is justified since the map
%ad identifies \CoOut$^-(A)$ with a subgroup of
%\CoOut$(A)$. 

We shall relate the groups ${\rm CoOut}^-(A)$ and $H_L^2(A)$.
In order to do that, we define an action of $\CoOut(A)$ on $H_L^2(A)$.
%Recall that any element $\mu \in \Reg^1(A)$ defines 
%an automorphism ad($\mu) = \mu^{-1}*{\rm id} * \mu$ 
%of the coalgebra $A$. In this way we get a group morphism
%${\rm ad} : \Reg^1(A) \longrightarrow {\rm Aut}_{\rm cog}(A)$.
%The group \CoInn$(A)$ of coinnner automorphisms
%of $A$ is then defined as 
%${\rm \CoInn}(A) = {\rm ad}({\rm Alg}(A,k))$. It is a normal subgroup
%of ${\rm Aut}_{\rm Hopf}(A)$ and \CoOut$(A)$ is defined
%to be the quotient group
%${\rm Aut}_{\rm Hopf}(A)/ {\rm \CoInn}(A)$.

Let $\sigma \in \Reg^2(A)$ and 
$\alpha \in \aut$. 
We put
$\sigma \gets \alpha := \sigma \circ (\alpha \otimes \alpha)$. 
It is clear that $\sigma \gets \alpha \in \Reg^2(A)$. 
The basic properties of this construction,
recorded in the following lemma,
enable us to construct a right action by automorphisms of 
$\CoOut(A)$ on $H_L^2(A)$.
 
\begin{lemm}\label{action}
Let $\alpha \in {\rm Aut}_{\rm Hopf}(A)$ and let
$\sigma \in \Reg^2(A)$.
\begin{enumerate}
\item If $\sigma \in \Reg_L^2(A)$, then $\sigma \gets \alpha \in 
\Reg_L^2(A)$.
\item If $\sigma$ is a left cocycle, then so is $\sigma \gets \alpha$.
\item Let $\omega \in \Reg^2(A)$. Then $(\sigma * \omega) \gets \alpha 
= 
(\sigma \gets \alpha)*(\omega \gets \alpha)$.
\item Let $\beta \in \aut$. Then
$\sigma \gets (\alpha\circ \beta) = (\sigma \gets \alpha)\gets \beta$.  
\item If $\mu \in \Reg_L^1(A)$, then $\mu \circ \alpha \in 
\Reg_L^1(A)$.
\item For $\mu \in \Reg^1(A)$, we have 
$\partial(\mu)\gets \alpha = \partial(\mu \circ \alpha)$.
\item If $\alpha\in {\rm \CoInn}(A)$ and $\sigma \in \Reg_L^2(A)$, then
$\sigma \gets \alpha$ = $\sigma$. 
\end{enumerate}
Hence the formula
$$
\begin{array}{rl}
H_L^2(A) \times {\rm CoOut}(A)& \longrightarrow  H_L^2(A) \\
(\overline{\sigma},\overline{\alpha})& \longmapsto   
\overline{\sigma \gets \alpha}
\end{array}
$$ 
defines a right action by automorphisms of ${\rm CoOut}(A)$
on $H_L^2(A)$.
\end{lemm}
\pf The proof of all statements follows by direct computation and we
leave it to the reader.\hfill$\Box$

\medskip

We can thus consider the semi-direct product group
${\rm CoOut}(A) \ltimes H_L^2(A)$.

\medskip

%We now construct another group attached to the Hopf algebra
%$A$, invoked in an exact sequence relating
%${\rm CoOut}(A) \ltimes H_L^2(A)$ and the biGalois group of $A$.
%Consider the subgroup
%${\rm ad}^{-1}({\rm Aut}_{\rm Hopf}(A))$ of $\Reg^1(A)$. Since 
%\CoInn$(A)$ is a normal subgroup of $A$, it follows that
%${\rm ad}^{-1}({\rm \CoInn}(A))$ is a normal subgroup of
%${\rm ad}^{-1}({\rm Aut}_{\rm Hopf}(A))$.
%
%
%
%\begin{lemm} Let $\mu \in \Reg^1(A)$. Then the following 
%assertions are equivalent.
%\begin{enumerate}
%
%\item {\rm ad}$(\mu) \in {\rm Aut}_{\rm Hopf}(A)$.
%
%
%\item {\rm ad}$(\mu^{-1}) \in {\rm Aut}_{\rm Hopf}(A)$.
%
%
%\item $\partial(\mu) \in \Reg_L^2(A)$.
%
%
%\item $\partial(\mu^{-1}) \in \Reg_L^2(A)$.
%
%\end{enumerate} 
\begin{lemma}\label{iota}The map  
$$
\begin{array}{rl}
i_0 : \cll &\longrightarrow 
{\rm CoOut}(A)\times H_L^2(A)\\
\mu &\longmapsto (\overline{{\rm
ad}(\mu)},\overline{\partial(\mu^{-1})})
\end{array}
$$ 
induces an injective group morphism
$$
\begin{array}{rl}
\iota :  {\rm CoOut}^-(A) 
&\longrightarrow {\rm CoOut}(A) \ltimes H_L^2(A)\\
 \overline{\mu} & \longmapsto 
(\overline{{\rm ad}(\mu)},\overline{\partial(\mu^{-1})}).
\end{array}
$$
\end{lemma}
\pf
%We have 
%$${\rm ad}(\mu) \circ m = 
%((\mu^{-1}) *{\rm id}_A * \mu) \circ m = 
%(\mu^{-1} \circ  m) * m * (\mu \circ m) \ {\rm and}$$ 
%$$m \circ ({\rm ad}(\mu) \otimes {\rm ad}(\mu)) =  
%m \circ ((\mu^{-1}* {\rm id}_A * \mu) \otimes (\mu^{-1}* {\rm id}_A * \mu))
%= (\mu^{-1} \otimes \mu^{-1}) * m * (\mu \otimes \mu).$$
%Hence we see that $(1) \iff (3)$ and $(2) \iff (4)$, while
%$(1) \iff (2)$ is obvious.
%Thus the map $i_0$ is defined.
%Since $Ker(\ad)=\Reg^1_L(A)$, the image of $\ad(\mu)$ does not depend
%on the choice of the cochain $\mu$, so $i_0$ is well-defined. 
The map $i_0$ is a group
morphism if 
$\partial(\lambda^{-1})\partial(\mu^{-1})=(\partial(\mu^{-1})\gets
\ad(\lambda))*\partial(\lambda^{-1})$ for every $\mu,\lambda\in\cll$.
By Lemma \ref{action},
$\partial(\mu^{-1}) \gets {\rm ad}(\lambda) = 
\partial( \mu^{-1}  \circ {\rm ad}(\lambda))$.
By Lemma \ref{list} and Lemma \ref{homom}, 
 $\partial(\mu^{-1} \circ {\rm ad}(\lambda)) * \partial(\lambda^{-1})
= \partial((\mu^{-1} \circ {\rm ad}(\lambda)) * \lambda^{-1})$. Using
that
$(\mu^{-1} \circ {\rm ad}(\lambda))* \lambda^{-1}
= \lambda^{-1} * \mu^{-1}$, we have that $i_0$ is a group morphism.

An element  $\gamma$ lies in the kernel of $i_0$ if and only if
$\gamma*\mu\in {\rm Alg}(A,k)$ for some
$\mu\in\Reg^1_L(A)$. Therefore, ${\rm Ker}(i_0)=\ad^{-1}(\CoInn(A))$ 
and the
statement is proved.\hfill$\Box$

\section{Examples}\label{exa}

This section contains, as first illustrative examples, the computation
of $H^2_L$ for Sweedler's Hopf algebra $H_4$ and for $E(n)$. The 
computation
here is based on the knowledge of explicit cocycles. These examples
will be generalized in the last sections using more theoretical
methods. Although $H_4$ coincides with $E(1)$, we deal with it
separately to have an easily manageable example throughout the paper.
In this Section we assume that $2$ is invertible in $k$.

\begin{example}\label{swee}{\rm
Let $H_4$ be Sweedler Hopf algebra,
with generators $g$  
and $x$, relations 
$$g^2=1,\quad x^2=0,\quad gx+xg=0
$$
and coproduct
$$\Delta(g)=g\otimes g,\quad \Delta(x)=1\otimes x+x\otimes g.
$$
It follows from self-duality that
%Then $\gamma\in \Reg^1(H_4)$ if  and  only if $\gamma(g)$ and 
%$\gamma(1)$ are
%invertible. Since the map $\Phi\colon H_4\to (H_4)^*$ with
%$\Phi(g)=1^*-g^*$ and $\Phi(h)=h^*+(gh)^*$ defines an Hopf algebra 
%isomorphism and since $Z(H_4)=k$, it follows that
$\Reg^1_{L}(H_4)$ is trivial. 

It is well known that 
${\rm Aut}_{\rm Hopf}(H_4)\cong k^\cdot$ where
$t\in k^\cdot$ corresponds to the automorphism $\alpha_t$ satisfying
$g\mapsto g$, $x\mapsto t
x$. Let $\gamma\in \Reg^1_{aL}(H_4)$. Then $\ad(\gamma)=\alpha_t$
for some $t\in k^\cdot$.
The relations $\ad(\gamma)(x)=tx$ and $\ad(\gamma)(gx)=tgx$ imply
that $\gamma(x)=\gamma(gx)=0$ and
$t=\gamma(g)\gamma(1)^{-1}=\gamma(1)\gamma^{-1}(g)=\pm1$. 
Hence, either $t=1$ and  $\gamma=\varepsilon$, or
$t=-1$ and $\gamma=(1^*-g^*)$. 
The map
$1^*-g^*$ is an algebra morphism 
so $\CoInn(H_4)=
\CoInt(H_4)\cong\regl1\cong
{\mathbb Z}_2$.
By Lemma \ref{equal} 
$$B^2_L(H_4)=\Reg^2_L(H_4)\cap 
\partial(\Reg^1(H_4))=\partial(\Reg^1_L(H_4))=\{\varepsilon^{\otimes 2}\}.$$ 

By the classification in \cite[Table]{[Ma]} the
lazy $2$-cocycles are exactly those appearing in
\cite{gio}, i.e., those $\sigma_t$ for $t\in k$ such that:
$$\sigma_t(g, g)=1,$$ 
$$\sigma_t(g, x)=\sigma_t(x, g)=\sigma_t(g,
gx)=\sigma_t(gx, g)=0$$
$$\sigma_t(x, x)=\sigma_t(gx, x)=-\sigma_t(x,
gx)=-\sigma_t(gx, gx)=\frac{t}{2}.$$
It is straightforward to check that $Z^2_{L}(H_4)\cong k$.
%The action of $Z^2_{L}(A)$ on the set of all crossed systems,
%parametrized as in \cite[Table]{[Ma]} by $t$-uples
%$(\alpha,\,\delta,\,u,\,a,\,b,\,s)$ is as follows. If $\sigma$
%corresponds to $(\alpha,\,\delta,\,u,\,a,\,b,\,s)$, then
%$\sigma*\sigma_{t}^{-1}$ will be the cocycle corresponding to
%$(\alpha,\,\delta,\,u,\,a-\frac{t}{2}1,\,b,\,s)$ so the orbits are
%parametrized by the $t$-uples $(\alpha,\,\delta,\,u,\,0,\,b,\,s)$.
%
%Finally,
Therefore
%$H^1_L(H_4)=\Reg_{L}^1(H_4)=1$ and 
$$H^2_L(H_4)\cong k.$$}
\end{example}

%
%\bigskip
\begin{example}\label{en}{\rm
Let $n$ be a positive integer and let $E(n)$ be the Hopf algebra 
generated
by $c$ and $x_i$ for $1\le i\le n$ with relations 
$$c^2=1;\quad cx_i+x_ic=0;\quad x_ix_j+x_jx_i=0;\quad x_i^2=0$$
and coproduct 
$$\Delta(c)=c\otimes c;\quad \Delta(x_i)=1\otimes x_i+x_i\otimes c.$$
%The 
%$1$-cochains are the elements $\gamma\in E(n)^*$ with
%$\gamma(c),\gamma(1)\in k^\cdot$. 

The map $\Phi\colon E(n)\to E(n)^*$
with $\Phi(1)=\varepsilon$, $\Phi(c)=1^*-c^*$,
$\Phi(x_i)=x_i^*+(cx_i)^*$ and $\Phi(cx_i)=x_i^*-(cx_i)^*$ for every
$i$ defines an Hopf algebra 
isomorphism. For a subset
$P=\{p_1,\ldots,p_l\}\subset\{1,\ldots,\,n\}$ we put $|P|=l$ and
$x_P=x_{p_1}\cdots x_{p_l}$. 
The centre of $E(n)$ is the span of $1$ and $x_P$ for $|P|$ even if
$n$ is odd and it is the span of $1$, $x_P$ for $|P|$ even and
$cx_1\cdots x_n$ if $n$ is even. 

It is not hard to check that $\Phi(x_P)=(x_P)^*+(cx_P)^*$ and that 
$\Phi(cx_1\cdots x_n)=(x_1\cdots
x_n)^*-(cx_1\cdots x_n)^*$. 
 Hence, $\Reg^1_L(E(n))$ consists of the
linear combinations of $\varepsilon$ (with coefficient $1$) and
elements of the above form. 

By \cite[Lemma 1]{PVO1} ${\rm Aut}_{\rm Hopf}(E(n))\cong {\rm 
GL}_n(k)$. The
automorphisms of $E(n)$ act trivially on $c$ and as linear maps on 
$x_1,\,\ldots,\,x_n$:
for $M\in {\rm GL}_n(k)$ the corresponding automorphism $\alpha_M$ is 
such
that $\alpha_M(c)=c$ and 
$\alpha_M(x_i)=\sum_j m_{ij}x_j$ for every $i=1,\ldots,n$. 
Let $\gamma\in \Reg^1_{aL}(E(n))$. Then  
$\ad(\gamma)=\alpha_M$ for some $M\in {\rm GL}_n(k)$.
The computation of $\ad(\gamma)(x_i)$,
$\ad(\gamma)(cx_i)$ and the relation
$\ad(\gamma)(cx_i)=(\ad(\gamma)(c))(\ad(\gamma)(x_i))$ imply that
$\gamma(cx_i)=\gamma(x_i)=0$ for every $i$ and that $M=\pm I_n$. 
If $M=I_n$, $\gamma\in \Reg^1_L(E(n))$; if $M=-I_n$, $\gamma\in
(1^*-c^*)*\Reg^1_L(E(n))$. Since $(1^*-c^*)$ is an algebra
morphism,  $\CoInn(E(n))=\CoInt(E(n))\cong{\mathbb Z}_2$ and
$Z^2(E(n))\cap \Reg^2_L(E(n))=B^2_{L}(E(n))$. 

Left 2-cocycles up to usual cohomology are classified in
  \cite{PVO2}. These classes are parametrized by $2$-cocycles
  satisfying some recurrence relations and
$$\sigma(c, c)=\alpha\in k^\cdot;\quad \sigma(c,
  x_i)=0;$$
$$\sigma(x_i, c)=\gamma_i;\quad\sigma(x_i, x_j)=m_{ij}$$
where $M_{ij}\in T_n(k)$, the vector space 
of lower-triangular $n\times n$ matrices with entries
in $k$.
It was proved in \cite[Lemma 2.2, Lemma 2.3]{GioJuan3}
% (up to a
%misprint in upper-lower triangular matrices), 
that these particular 
cocycles are lazy if and only if $\alpha=1$ and $\gamma_i=0$ for every
$i$ and that 
$\alpha=1$ and $\gamma_i=0$ for every $i$ is also a necessary condition 
for all
$2$-cocycles to be lazy. Let $\theta\in Z^2_L(E(n))$ and let 
$l_{ij}=\theta(x_i, x_j)$. If we consider the 
lazy cochain
$\gamma_\theta=\varepsilon+\sum_{i<j}l_{ij}((x_ix_j)^*+(cx_ix_j)^*)$ 
it is not hard to check that the lazy $2$-cocycle
$\theta*\partial(\gamma_\theta)(x_i, x_j)=0$ if $i<j$, while
 $\theta*\partial(\gamma_\theta)(x_i, x_j)=\theta(x_i,
x_j)+\theta(x_j, x_i)$ if $i>j$ and  
$\theta*\partial(\gamma_\theta)(x_i, x_i)=\theta(x_i, x_i)$. 
In particular, if $\theta$ is one of the cocycles constructed in
\cite{PVO2} and it is lazy, then $\gamma_\sigma=\varepsilon$.
The map 
$$
\begin{array}{rl}
\Psi\colon Z^2_L(E(n))&\longrightarrow T_{n}(k)\\ 
\sigma&\mapsto M_{ij}:=(\sigma*\partial(\gamma_{\sigma})(x_i, x_j))
\end{array}
$$
is surjective by \cite[Lemma 2.3]{GioJuan3}. 
Besides, if $\sigma,\omega\in Z_L^2(E(n))$ then 
$$\sigma*\partial(\gamma_\sigma)*\omega*\partial(\gamma_\omega)
=(\sigma*\omega)*\partial(\gamma_\sigma)*\partial(\gamma_\omega).
$$
%$$
%\begin{array}{l}
%\sigma*\partial\gamma_\sigma*\omega*\partial\gamma_\omega\\
%\phantom{\partial}=(\sigma*\omega)*\partial\gamma_\sigma*\partial\gamma_\omega.
%\end{array}
%$$
One has 
$$\sigma*\omega(x_i, x_j)=\sigma(x_i,
x_j)+\omega(x_i, x_j)$$ and 
%even if in general 
%$\gamma_\sigma*\gamma_\omega\neq \gamma_{\sigma*\omega}$ (it need not
%be zero on $x_1x_2x_3x_4$), 
%$$\partial(\gamma_\sigma*\gamma_\omega)(x_i, x_j)$$
$$\partial(\gamma_\sigma)*\partial(\gamma_\omega)(x_i,
x_j)=\partial(\gamma_\sigma)(x_i, x_j)+\partial(\gamma_\omega)(x_i, 
x_j)$$
and therefore $\Psi$ is a group morphism
$\Psi\colon Z^2_L(E(n))\to T_{n}(k)$. 

It is not hard to check that $\theta\in{\rm Ker}(\Psi)$ if and only if
$M=(\theta(x_i, x_j))$ is a skew-symmetric matrix, so that
that $\Psi(\partial(\gamma))=0$ for every $\partial(\gamma)\in 
B^2_L(E(n))$. 

On the other hand, if $\theta\in {\rm Ker}(\Psi)$ then the $2$-cocycle
$\omega:=\theta*\partial(\gamma_\theta)$
coincides with $\varepsilon\otimes\varepsilon$ when restricted to the
two-fold tensor product of the
span of $1$, $c$ and $x_1,\ldots,x_n$.  Hence, the cleft extension
$k\sharp_\omega E(n)$ is generated by $C$, $X_1,\ldots, X_n$ with
relations
$$C^2=1;\quad X_iX_j+X_jX_i=0;\quad CX_i+X_iC=0.$$
Since $k\sharp_\omega E(n)$  is isomorphic as an $E(n)$-module algebra 
to
$E(n)$, by \cite[Theorem 2.2]{doi} $\omega$ is a
(usual) coboundary, and it is a lazy $2$-cochain because $\theta$ and
$\partial(\gamma_\theta)$ are so. Hence, $\omega\in
B^2(E(n))\cap \Reg^2_L(E(n))=B^2_L(E(n))$ so
 $\theta\in B^2_L(E(n))$. Therefore
$\Psi$ defines a group isomorphism $H^2_L(E(n))\cong T_n(k)$. 
If we denote by $S^2(k^n)$ the group of $n\times n$ symmetric matrices
with entries in $k$ we have: 
$$H^2_L(E(n))\cong S^2(k^n).$$ 
This isomorphism will be
generalized in Section \ref{special} to a wider class of cotriangular Hopf 
algebras.}
\end{example}

\section{Lazy Galois objects}\label{galois}

It is classical to associate Hopf-Galois extensions to
2-cocycles. We characterize the 
Galois objects associated to lazy 2-cocycles as Galois objects
enjoying a special symmetry property: we will call them lazy Galois
objects.

Let us first recall a few facts concerning
Hopf-Galois extensions (\cite{[Mo]}). 
Let $A$ be a Hopf algebra.
A \textbf{right $A$-Galois extension (of $k$)}
is a non-zero right $A$-comodule algebra $Z$ with $Z^{{\rm co} A} = k$
such that the linear map
$\kappa_r$ defined by the composition
\begin{equation*}
\begin{CD}
\kappa_r : Z \otimes Z @>1_Z \otimes \rho>>
Z \otimes Z \otimes A @>m_Z \otimes 1_A>> Z \otimes A
\end{CD}
\end{equation*}
where $\rho$ is the coaction of $A$
and $m_Z$ is the multiplication of $Z$, is bijective.
We also say that a right $A$-Galois extension (of $k$)
is an \textbf{$A$-Galois object}.
A morphism of $A$-Galois objects is an $A$-colinear algebra 
morphism. 
It is known that any morphism of $A$-Galois extensions 
that are faithfully flat as $k$-modules is an isomorphism
(\cite[Remark 3.11]{[Sc]}).
The set of isomorphism classes of $k$-faithfully flat
$A$-Galois objects is denoted by Gal$(A)$. For example, 
when $A = k[G]$ is a group algebra, 
we have ${\rm Gal}(k[G]) \cong H^2(G,k^\cdot)$.
However ${\rm Gal}(A)$ does not carry a natural group structure in 
general. 

Galois objects are classically associated with 
2-cocycles in the following manner.
Let $\sigma$ be a left 2-cocycle. 
The right $A$-comodule algebra
$_{\sigma} \! A=k\sharp_\sigma A$ is defined in the following way. 
As a right $A$-comodule $_{\sigma} \! A = A$ and the product
of $_{\sigma}A$ is defined to be
$$a _{\sigma} \! . b = \sigma(a_{1}, b_{1}) a_{2} b_{2}, 
\quad a,b \in A.$$ 
The right $A$-Galois objects constructed from 2-cocycles are 
characterized
as the ones with the \textbf{normal basis property}, i.e.
those that are isomorphic to $A$ as $A$-comodules.
The Galois objects with the normal basis property are also 
characterized
as the \textbf{cleft} ones, see \cite{[Mo]} for details
and proofs of these statements. For future use, let us recall
how a left 2-cocycle is constructed from a right Galois
object having the normal basis property.
Let $Z$ be right $A$-Galois object with a right
$A$-colinear isomorphism $\psi : A \longrightarrow Z$ with $\psi(1)=1$.
Then $\sigma$ defined by $\sigma(a,b) =
\varepsilon(\psi ^{-1}(\psi(a)\psi(b)))$, $\forall a,b \in A$,
 is a left 2-cocycle
such that $\psi : {_{\sigma} \! A} \longrightarrow Z$
is an $A$-comodule algebra isomorphism.  

Let us now introduce an additional property
for Galois objects, that will be shown to be the
Galois translation of laziness for 2-cocycles.

\begin{defi}
A right $A$-Galois object $Z$ is said to be \textbf{lazy}
if there exists a right $A$-colinear isomorphism
$\psi : A \longrightarrow Z$ such that 
$\psi(1) = 1$ and such that the  morphism
$$\beta_\psi := (\psi^{-1} \otimes \psi) \circ \rho
: Z \longrightarrow A \otimes Z$$
is an algebra morphism. Such a map $\psi$ is called a \textbf{symmetry
morphism} 
for $Z$.
\end{defi}

The condition $\psi(1)=1$ is not a restriction. 
Given an isomorphism $\theta\colon Z\to M$ of right $A$-Galois
objects, $Z$ is lazy with symmetry morphism $\psi_Z$  if and only if 
$M$ is lazy with symmetry morphism $\psi_M=\theta\circ\psi_Z$. 
The subset of ${\rm Gal}(A)$ consisting of isomorphism classes
of lazy right $A$-Galois objects will be denoted by
${\rm Gal}^L(A)$. 

\begin{prop}\label{equivalence}
Let $Z$ be a right $A$-Galois object. Then the following
assertions are equivalent.
\begin{enumerate}
\item $Z$ is a lazy right $A$-Galois object.
\item There exists $\sigma \in Z_L^2(A)$ such that 
$_\sigma \! A \cong Z$ as right $A$-comodule algebras.
\end{enumerate}
\end{prop}
\pf $1 \Rightarrow 2$. Since $Z$ is lazy, it is cleft and we can assume
that $Z = {_{\omega} \!A}$ for a 2-cocycle $\omega$.
Let $\psi : A \longrightarrow {_\omega \! A}$ be a symmetry morphism.
Since $\psi$ is right $A$-colinear, there exists
$\mu=\varepsilon\circ\psi \in A^*$ such that 
$\psi(a) =  \mu(a_1)a_2$, $\forall a \in A$. The map $\mu\in\Reg^1(A)$
with inverse $\varepsilon\circ\psi^{-1}$. 
Then $\beta_\psi(a) = \mu^{-1}(a_1)a_2 \otimes \mu(a_3) a_4$.
Let $\sigma=\omega^\mu = (\mu \otimes \mu) * \omega *(\mu^{-1} \circ 
m)$.
By \cite[Theorem 2.2]{doi} $\sigma$ is a 2-cocycle
and $\psi$ induces a right $A$-comodule algebra
isomorphism ${_\sigma \! A} \longrightarrow {_\omega \! A}$.
Let us check that $\sigma$ is lazy.
For $a,b \in A$ we have
$$
\beta_\psi(a)\beta_\psi(b)
=\mu^{-1}(a_1)\mu^{-1}(b_1)\mu(a_3)\mu(b_3)\omega(a_4,b_4) a_2b_2 
\otimes a_5b_5 
$$
and 
$$\beta(a {_{\omega} \! \cdot}b) = \omega(a_1,b_1)
\mu^{-1}(a_2b_2) a_3b_3 \mu(a_4b_4)\otimes a_5b_5.$$   
Hence, since $\beta_\psi$ is an algebra morphism, we have
$$(\mu^{-1} \otimes \mu^{-1})*m*(\mu \otimes \mu)*\omega
= \omega*(\mu^{-1} \circ m)*m*(\mu \circ m).$$
This exactly means that
$m*\sigma = \sigma *m$ and thus $\sigma$ is lazy.

\noindent
$2\Rightarrow 1$. We can assume that $Z = {_\sigma \! A}$.
Taking $\psi = {\rm id}_A$ a direct computation
shows that $\beta_\psi = \Delta : {_\sigma \! A} \longrightarrow
A \otimes {_\sigma \! A}$ is an algebra morphism because
$\sigma$ is lazy.\hfill $\Box$

\medskip

\begin{rem}
{\rm It is possible to introduce a notion of lazy
Hopf-Galois extension for general $A$-Galois extensions
$R \subset Z$, corresponding to a notion of lazy
crossed system. Since this is not needed for the strict study of
the lazy cohomology group, it will done in the appendix.} 
\end{rem}

When $A$ is a $k$-flat Hopf algebra and
hence  $k$-faithfully flat, $\sigma\mapsto{_\sigma\!A}$ defines
 a surjective map $Z_L^2(A) \longrightarrow {\rm Gal}^L(A)$. 
By \cite[Theorem 2.2]{doi} this map induces a map
$H_L^2(A)$ to ${\rm Gal}^L(A)$, which is not
injective in general. This leads to consider biGalois objects. 

\medskip

Similarly to the right case,
a \textbf{left $A$-Galois object} 
is a non-zero left $A$-comodule algebra $Z$ with
$^{{\rm co}A}Z = k$ such that the linear map
$\kappa_l$ defined by the composition
\begin{equation*}
\begin{CD}
\kappa_l : Z \otimes Z @>\beta \otimes 1_Z>>
A \otimes Z \otimes Z @>1_A \otimes m_Z>> A \otimes Z
\end{CD}
\end{equation*}
where $\beta$ is the coaction of $A$ and $m_Z$ is the multiplication
of $Z$, is bijective.

Let $A$ and $B$ be Hopf algebras. An algebra $Z$ 
is said to be \textbf{an $A$-$B$-biGalois object} (cfr. \cite{[Sc1]}) 
if $Z$
is both a left $A$-Galois extension and a right $B$-Galois extension,
and if $Z$ is an $A$-$B$-bicomodule.
% When $A=B$, we simply say that
%$Z$ is an A-biGalois object. The set of isomorphism classes
%of $k$-faithfully flat $A$-$B$-BiGalois object is denoted by 
%BiGal$(A,B)$ (a morphism of biGalois objects being an 
%$A$-$B$-bicomodule algebras morphism),
%with BiGal$(A,A) = {\rm BiGal}(A)$.
 
Left Galois objects are related to 
2-cocycles in the following manner.
Let $\sigma$ be a right 2-cocycle. 
The left $A$-comodule algebra $A_\sigma$ is defined in the following 
way. 
As a left $A$-comodule $A_{\sigma} =A$ and the product
of $A_{\sigma}$ is defined to be
$$a \cdot_{\sigma} b = a_{1} b_{1}\sigma(a_{2}, b_{2}), 
\quad a,b \in A.$$
In the case of lazy 2-cocycles, we have
the following result. The simple proof is omitted.

\begin{prop}
Let $\sigma \in Z_L^2(A)$. Then ${_\sigma \! A}$,
endowed with $\Delta$ as left and right $A$-comodule structure, 
is a left $A$-Galois object, and is
an $A$-biGalois object. The corresponding biGalois object
is denoted by $A(\sigma)$.\hfill$\Box$ 
\end{prop}

More generally, we have the following result. 
%which is 
%a useful tool in explicit computations of 
%$H_L^2(A)$ once the set ${\rm Gal}(A)$ has been determined.
%
\begin{prop}
Let $Z$ be a lazy right $A$-Galois object. Then for any symmetry
morphism $\psi$, the map $\beta_\psi
: Z \longrightarrow A \otimes Z$ endows $Z$ with a left $A$-comodule
algebra structure for which $Z$ is $A$-biGalois.
\end{prop}
\pf We already know that $\beta_\psi$ is an algebra
morphism. Besides:
$$
\begin{array}{rl}
(\id_A\otimes\beta_\psi)\circ\beta_\psi&
=(\id_A\otimes\psi^{-1}\otimes\psi)\circ(\id_A\otimes\rho)\circ\beta_\psi\\
&=(\id_A\otimes\id_A\otimes\psi)\circ(\id_A\otimes\Delta)\circ(\id_A\otimes\psi^{-1})\circ\beta_\psi\\
&=(\id_A\otimes\id_A\otimes\psi)\circ(\id_A\otimes\Delta)\circ(\psi^{-1}\otimes\id_A)\circ\rho\\
&=(\id_A\otimes\id_A\otimes\psi)\circ(\id_A\otimes\Delta)\circ\Delta\circ\psi^{-1}\\
&=(\Delta\otimes\psi)\circ\Delta\circ\psi^{-1}\\
&=(\Delta\otimes\id_Z)\circ\beta_\psi
\end{array}
$$
and 
$$(\varepsilon\otimes\id_A)\circ(\psi^{-1}\otimes\psi)\circ\rho=(\varepsilon\otimes\psi)\circ\Delta\circ\psi^{-1}=\id_Z$$ 
hence
$\beta_\psi$ 
defines a left comodule structure on $Z$. Besides
$Z$ is an $A$-bicomodule because
$$
\begin{array}{rl}
(\beta_\psi\otimes\id_A)\circ\rho&=(\psi^{-1}\otimes\psi\otimes\id_A)\circ(\rho\otimes\id_A)\circ\rho\\
&=(\psi^{-1}\otimes\psi\otimes\id_A)\circ(\id_Z\otimes\Delta)\circ\rho\\
&=(\psi^{-1}\otimes\rho)\circ(\id_Z\otimes\psi)\circ\rho\\
&=(\id_A\otimes\rho)\circ\beta_\psi.
\end{array}
$$
Note also that $\beta_\psi \circ \psi = ({\rm id}_A \otimes \psi)
\circ \Delta$, hence $\psi$ is also a left $A$-colinear isomorphism.
There remains to check that $Z$ is left $A$-Galois for the coaction
$\beta_\psi$. 

If $f\colon Z\to T$ is a right $A$-comodule algebra isomorphism from
$Z$ to a lazy right $A$-Galois object $T$ for which the statement 
holds,
then $\varphi = f \circ \psi : A \longrightarrow T$
is a symmetry morphism and $\beta_\varphi$ endows $T$ with a left
comodule algebra morphism. A direct computation shows that
$f$ is also left $A$-colinear so $Z$ is
left $A$-Galois because $T$ is so.

Thus, by Proposition \ref{equivalence} 
we can assume that
$Z = {_\sigma \! A}$ for $\sigma \in Z_L^2(A)$.
Let $\psi : A\to {_\sigma \! A}$ be a symmetry 
morphism. As in the proof of Proposition \ref{equivalence} 
$\psi(a) = \mu(a_1)a_2$ for  every $a \in A$, with
$\mu=\varepsilon\circ\psi$ and
$\beta_\psi(a)=\mu^{-1}(a_1) a_2 \otimes\mu(a_3) a_4$.
Similarly to the proof of Proposition \ref{equivalence} we have 
$$(\mu^{-1} \otimes \mu^{-1})*m*(\mu \otimes \mu)*\sigma
= \sigma*(\mu^{-1} \circ m)*m*(\mu \circ m),$$
and since $\sigma \in Z_L^2(A)$, we see that
$\partial(\mu)*\sigma \in Z_L^2(A)$.
By standard theory
$\psi : {_{\partial(\mu) \sigma} \! A} \longrightarrow
{_\sigma \! A}$ is an algebra morphism and it is left colinear
(for $\Delta$ as a left coaction on ${_{\partial(\mu) \sigma} \! A}$).  
We conclude that ${_\sigma \! A}$ is left $A$-Galois for
the coaction $\beta_\psi$. \hfill $\Box$

\begin{defi}
Let $Z$ be an $A$-biGalois object. We say that $Z$ has \textbf{the 
binormal
basis property} if $Z \cong A$ as an $A$-bicomodule.
We also say that an $A$-biGalois object with
the binormal basis property is \textbf{bicleft}.
\end{defi}

We characterize now the biGalois objects arising 
from lazy 2-cocycles.

\begin{prop}\label{bicleft}
Let $Z$ be an $A$-biGalois object. Then the following
assertions are equivalent:
\begin{enumerate}
\item $Z$ has the binormal basis property.
\item  There exists $\sigma \in Z_L^2(A)$ such that 
$A(\sigma) \cong Z$ as $A$-bicomodule algebras.
\end{enumerate}
When this occurs, $Z$ is a lazy right $A$-Galois object.
\end{prop}
\pf $1 \Rightarrow 2$. Let $\psi : A \longrightarrow Z$ be an 
$A$-bicolinear isomorphism. We can assume without loss of generality
that $\psi(1) = 1$.
Let $\sigma : A \otimes A \longrightarrow k$ be defined
by $\sigma(a,b) = \varepsilon(\psi ^{-1}(\psi(a)\psi(b)))$.
Since $\psi$ is right $A$-colinear and $Z$ is right $A$-Galois,
it is well known that $\sigma$ is a left 2-cocycle and that
$\psi : {_\sigma \! A} \longrightarrow Z$ is a right $A$-comodule
algebra isomorphism. Similarly since
$\psi$ is left $A$-colinear and $Z$ is left $A$-Galois,
it is well known that $\sigma$ is a right 2-cocycle and that
$\psi : A_\sigma \longrightarrow Z$ is a left $A$-comodule algebra
isomorphism. Hence
$$\psi(a)\psi(b) = \psi(\sigma(a_1,b_1) a_2b_2) = 
\psi(\sigma(a_2,b_2)a_1b_1).$$
Since $\psi$ is bijective we
 conclude that $\sigma \in Z_L^2(A)$. It is then clear
that $\psi : A(\sigma) \to Z$ is an $A$-bicomodule
algebra isomorphism. 

\noindent
$2 \Rightarrow 1$. This follows because $A(\sigma) = A$ as an 
$A$-bicomodule. \hfill$\Box$

\medskip

Let us denote the set of isomorphism classes 
of $k$-faithfully flat $A$-$B$-biGalois objects by 
${\rm BiGal}(A,B)$ (a morphism being an $A$-$B$-bicomodule morphism),
and let ${\rm BiGal}(A) = {\rm BiGal}(A,A)$.
The isomorphism class of a faithfully flat $A$-biGalois object $Z$
is denoted by $[Z]$ in ${\rm BiGal}(A,B)$. 

We assume now that $A$ is a $k$-flat
Hopf algebra. Then BiGal$(A)$ inherits a natural group structure:
this is the \textbf{biGalois group} of $A$, defined
by Schauenburg \cite{[Sc1]}. The group law is induced by the
cotensor product: if $V$ is a right $A$-comodule and $W$ is a left
$A$-comodule 
their cotensor product over $A$, denoted
by $V\Box_AW$, is defined to be the kernel of
the linear map $\rho_V \otimes {\rm id}_W -  {\rm id}_V \otimes\rho_W :
V \otimes W \longrightarrow V \otimes A \otimes W$. In particular when
$V$ and $W$ are $A$-bicomodules, their cotensor product will be so. 
The faithful flatness assumption ensures
that we indeed have an associative law.
We denote by ${\rm Bicleft}(A)$ the subset of 
${\rm BiGal}(A)$ consisting of isomorphism classes of 
bicleft biGalois objects. We have the following result. 

\begin{prop}\label{normal}
Let $A$ be a $k$-flat Hopf algebra. 
Then {\rm Bicleft}$(A)$ is a normal subgroup
of {\rm BiGal}$(A)$.
\end{prop}
\pf The neutral element for the cotensor product $[A]$ lies in ${\rm
Bicleft}(A)$ which is stable for this group law. By Proposition
\ref{bicleft}, for 
any $[Z] \in {\rm Bicleft}(A)$, we have  $[Z] = [A(\sigma)]$ for
some lazy 2-cocycle $\sigma$.
Since $A(\sigma^{-1})\Box_A A(\sigma)\cong A$ we have 
$[A(\sigma^{-1})] = [A(\sigma)]^{-1}$ so
${\rm Bicleft}(A)$ is stable under inverses and it is a subgroup of 
${\rm
BiGal}(A)$. 
Now let $[Z] \in {\rm BiGal}(A)$ and $[T] \in {\rm Bicleft}(A)$.
We denote by  $Z^{-1}$ an $A$-biGalois object such
that $[Z^{-1}]= [Z]^{-1}$. Then we have, {\em as $A$-bicomodules}, 
$$Z  \square_A T \square_A Z^{-1} \cong
Z  \square_A A \square_A Z^{-1} \cong
Z  \square_A Z^{-1} \cong A$$
and hence Bicleft$(A)$ is normal in ${\rm BiGal}(A)$. \hfill$\Box$
  
\medskip

Most of the work has been done now to prove the following result. 

\begin{theo}
Let $A$ be a $k$-flat Hopf algebra. Then 
we have a group isomorphism
$$H_L^2(A) \cong {\rm Bicleft}(A).$$ 
In particular we may  identify $H_L^2(A)$ with
a normal subgroup of {\rm BiGal}$(A)$.
\end{theo}
\pf By Proposition \ref{bicleft} the map
$$
\begin{array}{rl}
f : Z_L^2(A) & \longrightarrow {\rm Bicleft}(A) \\
\sigma &\longmapsto [A(\sigma)]
\end{array}
$$
is surjective.
Using the coproduct $\Delta : A \longrightarrow A \otimes A$, 
we see by direct computation that for 
$\sigma$, $\omega \in Z_L^2(A)$ we have
$A(\sigma) \square_A A(\omega) \cong A(\sigma*\omega)$ and hence our 
map $f$ is a morphism of groups. A lazy $2$-cocycle $\sigma$ lies in 
${\rm Ker}(f)$ if and only if
there exists a biGalois isomorphism $\phi\colon A\to A(\sigma)$. Since 
$\phi$
is left and right colinear, we have
$\varepsilon(\phi(\a1))\a2=\phi(a)=\a1\varepsilon(\phi(\a2))$. Since
$\phi$ is an algebra morphism,
$\sigma=\partial(\varepsilon\circ\phi^{-1})$, hence ${\rm 
Ker}(f)\subset
B^2_L(A)$. Viceversa, if $\sigma=\partial(\gamma)$  with $\gamma\in
\Reg^1_L(A)$, then $a\mapsto\gamma(\a1)\a2$ gives an isomorphism
$A(\sigma)\to A$. \hfill$\Box$

\begin{coro}\label{a-b}
Let $A$ and $B$ be $k$-flat Hopf algebras. If there exists a 
$k$-faithfully 
flat $A$-$B$-biGalois object, then we have a group isomorphism
$$H_L^2(A) \cong H_L^2(B).$$
\end{coro}
\pf Let $Z$ be an 
$A$-$B$-biGalois object with $Z^{-1}$ its inverse in the 
Harrison groupoid $\mathcal H$ defined by Schauenburg in \cite[Section
4]{[Sc1]}, so that 
$Z^{-1}$ is a $B$-$A$-biGalois object. 
Using the composition law for biGalois objects
(in \cite[Section 4]{[Sc1]}) we have a map
\begin{align*}
{\rm BiGal}(A) & \longrightarrow {\rm BiGal}(B) \\
[T] & \longmapsto [Z^{-1} \Box_A T \Box_A Z]
\end{align*} 
that is clearly a group isomorphism. It is also clear that this map
induces an isomorphism between the groups of the corresponding bicleft 
objects
(same proof as in the proof of Proposition \ref{normal}), so we have 
our result by the previous theorem. \hfill$\Box$ 

\begin{rem}
{\rm The above result is a monoidal co-Morita invariance type result. 
Indeed
for two $k$-flat Hopf algebras, there exists a faithfully flat 
$A$-$B$ biGalois object if and only if
the linear monoidal categories  ${\rm Comod}(A)$ and ${\rm Comod}(B)$
are equivalent (\cite[Corollary 5.7]{[Sc1]}).}
\end{rem}

\begin{rem}\label{conjugate}{\rm
It follows from Corollary \ref{a-b} that the second  lazy cohomology
group is invariant under Doi's twists when $A$ is $k$-flat. This 
fact holds in general and it can be directly seen.
Indeed, if $\omega$ is a
$2$-cocyle for $A$, then $\sigma$
  is a lazy $2$-cocycle for the twisted Hopf algebra $_\omega
  A_{\omega^{-1}}$ if and only if $\omega^{-1}*\sigma*\omega$ is a
  lazy cocycle for $A$.   
%It is well known that if $\rho$ is a $2$-cocycle for $H$ and
%$\beta$ is a $2$-cocycle for $_\alpha H_{\alpha^{-1}}$ then
%$\alpha*\beta$ is a $2$-cocycle for $H$ and that, in particular,
%$\omega^{-1}$ is a $2$-cocycle for  $_\omega
%  H_{\omega^{-1}}$. It follows that  if $\sigma$
%  is a lazy $2$-cocycle for the twisted Hopf algebra $_\omega
%  H_{\omega^{-1}}$ then  $\omega^{-1}*\sigma*\omega$ is a
%  cocycle for $H$. Besides, it is clear that $\sigma$ centralizes
%  $\omega*m^*(H_L^*)*\omega^{-1}$ if and only if
%  $\omega^{-1}*\sigma*\omega$ centralizes $m^*(H^*)$, hence the
%  statement.
By centrality the assignment $\sigma\mapsto \omega^{-1}*\sigma*\omega$
  induces a group isomorphism $H^2_L(_\omega
  A_{\omega^{-1}})\to H^2_L(A)$.}
%
%By Lemma \ref{conjugation} the assignment defines a group
% isomorphism $$Z^2_L(_\omega
%  H_{\omega^{-1}},\f)\to Z^2_L(A).$$
% Besides, if $\sigma\in B^2_L(_\omega
%  H_{\omega^{-1}},\f)$ then $\sigma=D_\omega^1\gamma$
% for some $\gamma$ central in $(_\omega
%  H_{\omega^{-1}})^*$. Hence
%$$
%\begin{array}{rl}
%\sigma(h\otimes
%k)&=\sum\gamma(\h1)\gamma(\k1)\gamma^{-1}(\h2\cdot_\omega\k2)\\
%&=\sum\gamma(\h1)\gamma(\k1)\omega(\h2\otimes\k2)\gamma^{-1}(\h3\k3)\omega^{-1}(\h4\otimes\k4)
%\end{array}
%$$
%
%
%$\gamma\in H^*$($=(_\omega
%  H_{\omega^{-1}})^*$ as algebras) we get: 
%$\omega^{-1}*\sigma*\omega=D^1\gamma\in B^2_L(A)$. By symmetry, the
%  converse is also true, hence $H^2_L(A)\cong H^2_L(_\omega
%  H_{\omega^{-1}},\f)$.
%???? Write the explicit correspondence at the lazy cocycles level
%when the $A$-$B$ biGalois object is cleft, i.e. when
%$B$ is a twist of $A$ ???? 
%(still interesting, no flatness assumption)
\end{rem}

\begin{rem}\label{auto-equivalence}
{\rm Let us assume that $A$ is flat.
 According to the results of \cite{[Sc1]}, the group
${\rm BiGal}(A)$ is isomorphic to the group of 
isomorphism classes of $k$-linear monoidal auto-equivalences of 
${\rm Comod}(A)$. 
The subgroup ${\rm Bicleft}(A)$ is then, in this setting, identified
with the group of isomorphism classes of $k$-linear monoidal
 auto-equivalences of  
${\rm Comod}(A)$ that are isomorphic, as linear functors, with the 
identity
functor. This might be checked directly at the lazy cocycle level.
With this interpretation of lazy cohomology,  
since the category of comodules
over a coquasibialgebra is a monoidal category, we recover
a natural way to define the lazy cohomology groups for
coquasibialgebras as done by Schauenburg in \cite{[Sc3]}.}   
\end{rem}

For a cocommutative ($k$-flat) Hopf algebra $A$, there is a group 
isomorphism
${\rm Aut}_{\rm Hopf}(A) \ltimes {\rm Gal}(A) \cong {\rm BiGal}(A)$
\cite[Section 5]{[Sc1]}. A group morphism between these two groups
 still exists in general and 
it induces the following exact sequence. 

\begin{theo}
Let $A$ be a $k$-flat Hopf algebra. There is a group exact sequence
$$
1 \longrightarrow {\rm CoOut}^-(A)\longrightarrow 
{\rm CoOut}(A) \ltimes H_L^2(A) \longrightarrow {\rm BiGal}(A).
$$
\end{theo}
\pf To any pair $(\alpha,\sigma) \in {\rm Aut}_{\rm Hopf}(A) \times 
Z^2_L(A)$,
we associate the $A$-biGalois object ${^\alpha \! A}(\sigma)$, which
is $A(\sigma)$ as a right $A$-comodule algebra with
left coaction $\rho = (\alpha \otimes {\rm id}_A) \circ \Delta$.   
This assignment induces a well-defined map
\begin{align*}
\Upsilon : {\rm CoOut}(A) \ltimes & H_L^2(A)  \longrightarrow {\rm 
BiGal}(A) \\
(\overline{\alpha}, & \overline{\sigma})  \longmapsto [{^\alpha 
\!A}(\sigma)].
\end{align*}
Indeed, if $\mu \in {\rm Reg}^1_L(A)$ and $\phi \in {\rm Alg}(A,k)$, 
the map $\phi^{-1} * \mu * {\rm id}_A$ is an isomorphism
of biGalois objects between the $A$-bicomodule algebras  
${^{\alpha \circ {\rm ad}(\phi)} \! A}(\sigma \partial(\mu))$ and
${^\alpha \! A}(\sigma)$.  
%Hence we have a well-defined map
%\begin{align*}
%\Upsilon : {\rm CoOut}(A) \ltimes & H^2(A)  \longrightarrow {\rm BiGal}(A) \\
%(\overline{u}, & \overline{\sigma})  \longmapsto [{^u \!A}(\sigma)]
%\end{align*}

Let $(\alpha,\sigma), (\beta,\omega) \in {\rm Aut}_{\rm Hopf}(A) \times 
Z^2_L(A)$. 
Let $\gamma : A \rightarrow A \otimes A$ be defined by
$\gamma(a) = \beta(a_1) \otimes a_2$. It is straightforward to check
that $\gamma$ induces an $A$-bicomodule algebra map
${^{\alpha\circ \beta} \! A}((\sigma \gets \beta)* \omega) \rightarrow
{^\alpha \! A}(\sigma) \square_A {^\beta \! A}(\omega)$ and hence an 
isomorphism.
Thus $\Upsilon$ is a morphism of groups.
Let $(\overline{\alpha},\overline{\sigma}) \in {\rm Ker}(\Upsilon)$:
there exists an $A$-bicomodule algebra isomorphism
$f : {^\alpha \! A}(\sigma) \rightarrow A$. Since $f$ is right 
$A$-colinear,
there exists $\mu \in {\rm Reg}^1(A)$ such that $f =\mu * {\rm id}_A$.
Since $f$ is left colinear, we have $\mu * {\rm id}_A = \alpha * \mu$
and hence $\alpha = {\rm ad}(\mu^{-1})$.
Since $f$ is an algebra map, we have $\partial(\mu)= \sigma \in 
Z^2_L(A)$.
Therefore $(\overline{\alpha},\overline{\sigma}) =
\iota(\overline{\mu^{-1}})$ where $\iota$ is as in 
Lemma \ref{iota}. Therefore ${\rm Ker}(\Upsilon)\subset \iota({\rm
CoOut}^-(A))$. Viceversa, if $\theta\in \cll$, then
the map $\theta^{-1}*\id_A$ is a bicomodule algebra isomorphism
$^{\ad(\theta)}\!A(\partial(\theta^{-1}))\to A$, whence the statement. 
\hfill$\Box$

\medskip

\begin{example}{\rm Let $A$ be Sweedler's Hopf algebra $H_4$. 
By the results in Example \ref{en}, ${\rm Ker}(\iota)$ is trivial and 
 ${\rm CoOut}(H_4)\ltimes H^2_L(H_4)\cong k^\cdot/(\pm1)\ltimes 
  k$.
%with product $(a,t)(b,r)=(ab, tb^2r)$. 
%because  $\sigma_t\gets\alpha_b=\sigma_{tb^2}$. 
By \cite[Theorem 5]{[Sc2]} ${\rm BiGal}(H_4)\cong
  k^\cdot\ltimes k$ so the map $\Upsilon$ is not surjective and
 ${\rm CoOut}(H_4)\ltimes H^2_L(H_4)$ is a normal subgroup of
 ${\rm BiGal}(H_4)$. Their quotient is $k^\cdot/(k^\cdot)^2\cong 
{H}^2({\mathbb Z}_2, k^\cdot)$.  
The subgroup
  $k$ corresponding to $(1,\,t)$ in \cite[Theorem 5]{[Sc2]}
  corresponds to $H^2_L(H_4)$ while 
  the subgroup ${\rm CoOut}(H_4)$ corresponds to the elements in 
${\rm BiGal}(H_4)$
  of the form $(a^2, 0)$ for $a\in k^\cdot$.
Let us observe that ${\rm CoOut}(A)\ltimes H^2_L(A)$ 
appeared in \cite{embedding} as a subgroup of the Brauer group 
$BQ(k, H_4)$.} 
\end{example}

\section{The Kac-Schauenburg exact sequence}\label{KS}

This section is devoted to the construction of a Kac-Schauenburg-type 
exact sequence for a double crossed Hopf algebra of
possibly infinite-dimensional Hopf algebras.
As a consequence we will derive a Schur-Yamazaki type formula,
as well as the description of the lazy cohomology
of the Drinfeld double of (finite-dimensional) Hopf algebras.

\medskip

Recall \cite{[KS]}
that a matched pair $(B,A, \rightharpoonup, \leftharpoonup)$ consists
of two Hopf algebras $B$ and $A$ together with linear maps
$\rightharpoonup  : A \otimes B \longrightarrow B$ and
$\leftharpoonup :  A \otimes B \longrightarrow A$ making 
$B$ into a left $A$-module coalgebra and $A$ into a right $B$-module coalgebra
respectively, and satisfying the following conditions

$$ (aa')\leftharpoonup b = (a \leftharpoonup(a'_1 \rightharpoonup b_1)) 
(a'_2 \leftharpoonup b_2),
\quad 1\leftharpoonup b = \varepsilon_B(b) 1,$$
$$a \rightharpoonup (bb') = (a_1 \rightharpoonup b_1)
((a_2 \leftharpoonup b_2)\rightharpoonup b'), \quad
a\rightharpoonup 1 = \varepsilon_A(a)1 , $$
$$ (a_1 \leftharpoonup b_1) \otimes (a_2 \rightharpoonup b_2)
= (a_2 \leftharpoonup b_2) \otimes (a_1 \rightharpoonup b_1). $$

One associates a Hopf algebra $\bic$
to a matched pair $(B,A, \rightharpoonup, \leftharpoonup)$ in the
following manner: as a coalgebra $\bic$ is the tensor product
coalgebra $B\otimes A$, the product (with unit
$1 \otimes 1$) is defined by
$$(b \otimes a) (b' \otimes a')
= b (a_1 \rightharpoonup b'_1) \otimes (a_2 \leftharpoonup b'_2) a',$$
and the antipode is defined by
$$S(a \otimes b) = (1 \otimes S_B(b)) (S_A(a)\otimes 1).$$
The Hopf algebra $\bic$ is said to be a double crossed product Hopf algebra.

\medskip

\begin{defi}
Let $(B,A, \rightharpoonup, \leftharpoonup)$ be a matched pair and 
consider the associated Hopf algebra $\bic$. A \textbf{central
$\bic$-pairing} is a convolution invertible linear map 
$\beta : B \otimes A \longrightarrow k$ satisfying the following
conditions: 
$$\beta(bb',a) = \beta(b_1,a_1) \beta(b',a_2 \leftharpoonup b_2) \quad , \quad
\beta(b,aa') = \beta(b_1,a'_1) \beta(a'_2\rightharpoonup b_2 , a)$$
$$\beta(b_2,a_2)a_1 \leftharpoonup b_1 = 
\beta(b_1,a_1) a_2 \leftharpoonup b_2 \quad , \quad
\beta(b_1,a_1) a_2\rightharpoonup b_2 = \beta(b_2,a_2) a_1 \rightharpoonup b_1$$
$$\beta(1,a) = \varepsilon_A(a) \quad , \quad 
\beta(1,b) = \varepsilon_B(b).$$
The set of central $\bic$-pairings is denoted by $\cp(\bic)$.
\end{defi}
 
It is straightforward to check that $\cp(\bic)$, endowed with the convolution
product, is a group. 

\medskip

Here is our version of the
Kac-Schauenburg exact sequence. A brief comparison with
\cite[Theorem 6.5.1]{[Sc3]} will be given after the proof. 

\begin{theo}\label{sequence} Let $(B,A, \rightharpoonup, \leftharpoonup)$ be a matched pair
of Hopf algebras and consider the associated 
double crossed product Hopf algebra
$\bic$. Then we have a group exact sequence:
\begin{equation*}
\begin{CD}
1 @>>> H^1_L(\bic) @> {\rm res}>> H^1_L(B) \times H^1_L(A)
@> \Lambda >> \cp(\bic) \\ @> \Sigma >>  H^2_L(\bic) @> {\rm res} >>
H^2_L(B) \times H^2_L(A)  
\end{CD}
\end{equation*}
\end{theo}

Since the Hopf algebras $A$ and $B$ are identified with Hopf subalgebras of $\bic$, we get well-defined restriction maps, which clearly are group
morphisms. Also it is clear that the first one is injective. We now have to construct
the maps $\Lambda$ and $\Sigma$.

\begin{lemm}\label{primo}
Let $(\phi_B,\phi_A) \in H^1_L(B) \times H^1_L(A)$. Define
$\Lambda(\phi_A,\phi_B) : B \otimes A \longrightarrow k$ by
$$\Lambda(\phi_B,\phi_A)(b,a) = \phi_A^{-1}(a_1) \phi_B^{-1}(b_1)
\phi_B(a_2 \rightharpoonup b_2) \phi_A(a_3\leftharpoonup b_3).$$
Then $\Lambda(\phi_A,\phi_B) \in \cp(\bic)$, and this defines a group morphism
such that the sequence
\begin{equation*}
\begin{CD}
H^1_L(\bic) @> {\rm res}>> H^1_L(B) \times H^1_L(A)
@> \Lambda >> \cp(\bic)
\end{CD} 
\end{equation*}
is exact.
\end{lemm}
\pf It is not difficult to check, using 
the conditions of a matched pair and the laziness of 
$\phi_B$ and $\phi_A$, that $\Lambda(\phi_B,\phi_A) \in \cp(\bic)$
and that $\Lambda$ is a group morphism.
For $\phi \in{\rm Reg}_L^1(\bic)$, we have
$\phi \in H^1_L(\bic)$ if and only if
the restrictions to $A$ and $B$
belong to $H^1_L(A)$  and $H^1_L(B)$, respectively, and 
$$\phi(b \otimes a) = \phi(b  \otimes 1) \phi(1 \otimes a)
=\phi((1\otimes a)(b\otimes 1))$$
$$=\phi(a_1 \lu b_1 \otimes a_2 \ru b_2) 
=\phi(a_1 \lu b_1 \otimes 1) \phi(1 \otimes a_2 \ru b_2).$$
Hence we see that for $\phi  \in H^1_L(\bic)$, 
then $\Lambda({\rm res}(\phi)) = \varepsilon_B \otimes \varepsilon_A$.
Conversely, if $(\phi_B,\phi_A) \in {\rm Ker}(\Lambda)$ then 
$\phi = \phi_B \otimes \phi_A \in H^1_L(\bic)$ and 
$(\phi_B,\phi_A) = {\rm res}(\phi)$. \hfill$\Box$ 

\bigskip

\begin{lemm} Let $\beta \in \cp(\bic)$. Let 
$\sigma_\beta : B \otimes A \otimes B \otimes A\longrightarrow k$ be defined
by 
$$\sigma_\beta(b \otimes a, b' \otimes a') = \beta(b',a)
\varepsilon_B(b) \varepsilon_A(a).$$
Then $\sigma_\beta \in Z^2_L(\bic)$ defining a group morphism
\begin{align*} 
\Sigma : \cp(\bic) & \longrightarrow H^2_L(\bic) \\
\beta & \longmapsto \overline{\sigma_\beta}
\end{align*}
\end{lemm}
\pf It is easy to see, using that $\beta \in \cp(\bic)$, that
$\sigma_\beta \in {\rm Reg}_L^2(\bic)$. Let us check that
$\sigma_\beta$ is a (left) 2-cocycle. We have
$$
\begin{array}{l}
\sigma_\beta(b_1 \otimes a_1, b'_1 \otimes a'_1)
\sigma_\beta((b_2 \otimes a_2)(b'_2\otimes a'_2), b'' \otimes a'')\\
\phantom{\sigma_\beta(b_1 \otimes a_1)}=\sigma_\beta(b_1 \otimes a_1, b'_1 \otimes a'_1)
\sigma_\beta((b_2 [a_2 \rightharpoonup b'_2] \otimes [a_3 \leftharpoonup b'_3]a'_2, 
b'' \otimes a'')\\
\phantom{\sigma_\beta(b_1 \otimes a_1)}=\varepsilon_B(b) \varepsilon_A(a'') \beta(b'_1,a_1) 
\beta(b'', [a_2 \leftharpoonup b'_2]a')\\
\phantom{\sigma_\beta(b_1 \otimes a_1)}=\varepsilon_B(b) \varepsilon_A(a'') \beta(b'_1,a_1) \beta(b''_1,a'_1) 
\beta(a'_2 \rightharpoonup b''_2, a_2 \leftharpoonup b'_2) \\ 
\phantom{\sigma_\beta(b_1 \otimes a_1)}=\varepsilon_B(b) \varepsilon_A(a'') \beta(b''_1,a'_1) 
\beta(b'[a'_2 \rightharpoonup b''_2], a) \\ 
\phantom{\sigma_\beta(b_1 \otimes a_1)}=\sigma_\beta(b'_1 \otimes a'_1, b''_1 \otimes a''_1)
\sigma_\beta(b \otimes a, b'_2 [a'_2 \rightharpoonup b''_2] \otimes 
[a'_3 \leftharpoonup b''_3]a''_2)\\
\phantom{\sigma_\beta(b_1 \otimes a_1)}=\sigma_\beta(b'_1 \otimes a'_1, b''_1 \otimes a''_1)
\sigma_\beta(b \otimes a,
(b'_2 \otimes a'_2) (b''_2\otimes a''_2)),
\end{array}
$$
where the properties of the central $\bic$-pairing $\beta$
have been used. Hence $\sigma_\beta \in Z^2_L(\bic)$ and it is immediate
that the induced map $\Sigma$ is a group morphism. \hfill$\Box$

\begin{lemm}\label{secondo} The sequence
\begin{equation*}
\begin{CD}
H^1_L(B) \times H^1_L(A)
@> \Lambda >> \cp(\bic)  @> \Sigma >>  H^2_L(\bic) 
\end{CD}
\end{equation*}
is exact.
\end{lemm}
\pf Let $(\phi_B,\phi_A) \in H^1_L(B) \times H^1_L(A)$, and let 
$\mu = \phi_B^{-1} \otimes \phi_A^{-1}$. Clearly $\mu \in {\rm Reg}_L^1(\bic)$,
and a direct computation shows that 
$\Sigma(\Lambda(\phi_B,\phi_A)) = \overline{\partial(\mu)}$.
Hence ${\rm Im}(\Lambda) \subset {\rm Ker}(\Sigma)$. 

Conversely, let $\beta \in {\rm Ker}(\Sigma)$. Then there exists 
$\mu \in {\rm Reg}^1_L(\bic)$ such that $\sigma_\beta = \partial(\mu)$.
Let $\phi_B : B \rightarrow k$ be defined by $\phi_B(b) = \mu^{-1}(b \otimes 1)$,
and similarly, let $\phi_A : A\rightarrow k$ be defined 
by $\phi_A(a) = \mu^{-1}(1 \otimes a)$. 
Since $\mu \in {\rm Reg}^1_L(\bic)$, 
we have $\phi_B \in {\rm Reg}^1_L(B)$ and $\phi_A \in {\rm Reg}^1_L(A)$.
Computing $\sigma_\beta(b \otimes 1,b' \otimes 1)$ and 
$\sigma_\beta(1 \otimes a,1 \otimes a')$, we see that
$\phi_B \in H^1_L(B)$ and $\phi_A \in H^1_L(A)$.
Computing $\sigma_\beta(b \otimes 1, 1 \otimes a)$, we find
$\mu^{-1}(b \otimes a) = \mu^{-1}(b \otimes 1) \mu^{-1}(1 \otimes a)$.
Then we have
\begin{align*}
\beta(b,a) & = \sigma_\beta(1 \otimes a , b \otimes 1)
=\mu(1 \otimes a_1) \mu(b_1 \otimes 1) 
\mu^{-1}(a_2 \rightharpoonup b_2 \otimes a_3\leftharpoonup b_3) \\
& = \phi_A^{-1}(a_1) \phi_B^{-1}(b_1)
\phi_B(a_2 \rightharpoonup b_2) \phi_A(a_3\leftharpoonup b_3),
\end{align*}
which proves that $\beta = \Lambda(\phi_B,\phi_A)$ and finishes the
proof of the lemma.\hfill$\Box$  

\bigskip

There remains to check the exactness of our sequence at $H_L^2(\bic)$. 
The following lemma, taken from \cite{[Sc3]}, will be useful:

\begin{lemm}\label{class}(\cite[Lemma 6.2.7]{[Sc3]})
Let $\sigma \in Z^2_L(\bic)$. Then there exists $\sigma' \in Z^2_L(\bic)$ having the
same class as $\sigma$ in $H^2_L(\bic)$ and satisfying:
$$\sigma'(b \otimes 1, b' \otimes a') =  
\sigma(b \otimes 1,b'\otimes 1) \varepsilon_A(a'), \ 
{\rm and} \
\sigma'(b \otimes a,1\otimes a') = \sigma(1 \otimes a,1 \otimes a')
\varepsilon_B(b).$$   
\end{lemm}
\pf Let $\mu : B \otimes A \longrightarrow k$ be defined 
by $\mu(b \otimes a) = \sigma(b \otimes 1, 1 \otimes a)$. It is straightforward
to check that $\mu \in {\rm Reg}^1_L(\bic)$ and that
$\sigma' = \sigma * \partial(\mu)$ satisfies the above
conditions. \hfill$\Box$ 

\begin{lemm}\label{terzo}
The sequence 
\begin{equation*}
\begin{CD}
 \cp(\bic)  @> \Sigma >>  H^2_L(\bic) @> {\rm res} >>
H^2_L(B) \times H^2_L(A) 
\end{CD}
\end{equation*}
is exact.
\end{lemm}
\pf It is clear from the definitions that
${\rm Im}(\Sigma) \subset {\rm Ker} ({\rm res})$. 
Conversely let $\sigma \in Z^2_L(\bic)$ be such that 
$\overline{\sigma} \in {\rm Ker}({\rm res})$. 
Then there exists $\mu_1 \in {\rm Reg}_L^1(B)$ and
 $\mu_2 \in {\rm Reg}_L^1(A)$ such that
$$\sigma(b\otimes 1, b'\otimes 1)=
\mu_1(b_1) \mu_1(b'_1) \mu_1^{-1}(b_2b'_2) \ {\rm and} \   
\sigma(1\otimes a, 1 \otimes a')=
\mu_2(a_1) \mu_2(a'_1) \mu_2^{-1}(a_2a'_2).$$
Define $\mu : B \otimes A \longrightarrow k$
by $\mu(b\otimes a) = \mu_1(b) \mu_2(a)$. Clearly
$\mu \in {\rm Reg}_L^1(\bic)$. A direct computation gives
$$\sigma*\partial(\mu^{-1})(b\otimes 1, b'\otimes 1)=
\varepsilon_B(bb') \quad {\rm and} \quad   
\sigma*\partial(\mu^{-1})(1\otimes a, 1 \otimes a')=
\varepsilon_A(aa').$$
Therefore we can assume, without changing
the class of $\sigma$ in $H^2_L(\bic)$, that
$$\sigma(b\otimes 1, b'\otimes 1)=
\varepsilon_B(bb') \quad {\rm and} \quad   
\sigma(1\otimes a, 1 \otimes a')=
\varepsilon_A(aa').$$
By Lemma \ref{class} we can assume that
$$\sigma(b\otimes 1, b'\otimes a')=
\varepsilon_B(bb') \varepsilon_A(a') \quad {\rm and} \quad   
\sigma(b\otimes a, 1 \otimes a')=
\varepsilon_B(b) \varepsilon_A(aa').$$
Then using the fact that $\sigma$ is a 2-cocycle, we have
$$
\begin{array}{rl}
\sigma(b\otimes a, b'\otimes a') & =   
\varepsilon_B(b_1) \varepsilon_A(a_1) \sigma(b_2\otimes a_2, b'\otimes a') \\
& = \sigma(b_1 \otimes 1, 1 \otimes a_1)
 \sigma(b_2\otimes a_2, b'\otimes a') \\
& = \sigma(1\otimes a_1, b'_1\otimes a'_1) 
\sigma(b\otimes 1, (1  \otimes a_2)(b'_2\otimes a'_2)) \\
& = \varepsilon_B(b) \sigma(1\otimes a, b'\otimes a')\\
& = \varepsilon_B(b) \sigma(1\otimes a, b'_2\otimes a'_2)
\sigma(b'_1\otimes 1, 1\otimes a'_1) \\
& = \varepsilon_B(b) \sigma(1\otimes a_1, b'_1\otimes 1)
\sigma((1 \otimes a_2)(b'_2 \otimes 1), 1\otimes a') \\
& = \varepsilon_B(b) \varepsilon_A(a') \sigma(1\otimes a, b'\otimes 1).
\end{array}
$$ 
Thus there remains to check that
$\beta : B \otimes A \longrightarrow k$ defined by
$\beta(b,a) = \sigma(1 \otimes a, b \otimes 1)$
belongs to $\cp(\bic)$. It is clear that $\beta$ is convolution invertible.
We have
\begin{align*}
 \beta(bb',a)& = \sigma(1\otimes a, bb'\otimes 1)
= \sigma(1\otimes a, (b_1\otimes 1)(b'_1 \otimes 1))
\sigma(b_2\otimes 1, b'_2\otimes 1) \\
& = \sigma(1\otimes a_1, b_1\otimes 1) 
\sigma(a_2 \lu b_2\otimes a_3 \ru b_3, b' \otimes 1)
= \beta(b_1,a_1) \beta(b',a_2 \ru b_2).
\end{align*}
Using that $\sigma$ is lazy for $\bic$ and the counits, 
 one sees easily that $\beta$ satisfies
the last two conditions defining central $\bic$-pairings. 
Finally one checks that
$\beta(b,aa') = \beta(b_1,a'_1) \beta(a'_2 \lu b_2,a)$
similarly to the case of the first condition, using
centrality of $\beta$.\hfill$\Box$ 

\bigskip

\noindent Combining Lemma \ref{primo}, \ref{secondo} and \ref{terzo} together concludes the proof of 
Theorem \ref{sequence}.\hfill$\Box$

\bigskip

Let us briefly compare 
Theorem \ref{sequence} with Schauenburg's version of G.I. Kac's exact sequence
\cite[Theorem 6.5.1]{[Sc3]}. For this we assume that
$A$ is finite-dimensional. Then the group
$\cp(\bic)$ is isomorphic to the automorphism group 
${\rm Aut}_{\rm ext}(A^* \# B)$ of the Hopf algebra extension
corresponding to the bismash product of the dual Singer pair of
the original matched pair. Therefore in this case 
we have the same exact sequence.

\bigskip

We now examine the case of a matched pair 
$(B,A, \rightharpoonup, \leftharpoonup)$ with 
$\lu$ and $\ru$ trivial. In this case the double crossed product
$\bic$ is the tensor product Hopf algebra $B \otimes A$.
In this case the group $\cp(B\otimes A)$ is clearly an abelian
group.
Examining the Kac-Schauenburg exact sequence, we get the following
Schur-Yamazaki type formula, generalizing the classical
one in group cohomology:

\begin{theo}\label{Yamazaki}
Let $A$ and $B$ be Hopf algebras.
Then we have a group isomorphism
$$H_L^2(A\otimes B) \cong H_L^2(A) \times H_L^2(B)
\times \mathcal{ZP}(A \otimes B).$$
\end{theo}
\pf
First note that since $\lu$ and $\ru$ are trivial,
then the map $\Lambda$ is trivial, and by the exactness at $\cp(B \otimes A)$
we see that $\Sigma$ is injective. Now
for $(\sigma_1, \sigma_2) \in  Z_L^2(A) \times Z_L^2(B)$, it 
is straigtforward to check that $\sigma : B \otimes A \longrightarrow k$
defined by $\sigma(b \otimes a,b'\otimes a') = \sigma_1(b,b') \sigma_2(a,a')$
is a lazy 2-cocycle. This defines a group morphism
$j : H_L^2(A) \times H_L^2(B) \longrightarrow H_L^2(B \otimes A)$
such that ${\rm res} \circ j = {\rm id}$.
Thus by Theorem \ref{sequence}
we have a split exact sequence: 
\begin{equation*}
\begin{CD}
1 @>>> \cp(B \otimes A)  @> \Sigma >>  H^2_L(B \otimes A) @> {\rm res} >>
H^2_L(B) \times H^2_L(A) @>>> 1  
\end{CD}
\end{equation*}
There just remains to be remarked that, since 
both actions are trivial,
$\Sigma$ maps $\cp(B \otimes A)$ into a central 
subgroup of $H^2_L(B \otimes A)$: our split exact sequence
is central and hence we have the announced
isomorphism. \hfill$\Box$ 

\medskip

We assume for the rest of the section that $k$ is a field.
We shall give an application of the previous formula to the computation
of the second lazy cohomology group of a Drinfeld double.

First let us begin with a more convenient description of the
group of central $A\otimes B$-pairings.
Let $A$ and $B$ be Hopf algebras. The set of Hopf algebra
morphisms $f : A \to B$ satisfying
\begin{equation}\label{ell}
f(A) \subset Z(B) \quad {\rm and} \quad
f(a_1)\otimes a_2 = f(a_2) \otimes a_1, \ \forall a \in A,
\end{equation}
will be denoted by $\mathcal L(A,B)$.
% (NAME FOR SUCH MORPHISMS ???). 
A direct computation shows that ${\mathcal L}(A,B)$ is a group under 
the convolution product (the inverse $f\circ S$ of $f$ is a Hopf 
algebra
morphism because of (\ref{ell})). We put ${\mathcal L}(A) :={\mathcal 
L}(A,A)$.

\begin{lemm}\label{pairings}
Let $A$ and $B$ be Hopf algebras. 
Assume that $B^o$, the dual Hopf algebra of $B$,
separates the points of $B$. Then the groups
$\mathcal{ZP}(A\otimes B)$ and $\mathcal L(A,B^o)$ are isomorphic.
\end{lemm}
\pf For any pairing $\beta : A \otimes B \rightarrow k$
we define the Hopf algebra map $f_\beta : A \rightarrow B^o$ by
$f_\beta(a)(b) := \beta(a,b)$. 
It is well known that this establishes a bijective correspondence
between pairings $A\otimes B \rightarrow k$ and Hopf algebra maps
$A \rightarrow B^o$. 
Let $a \in A$ and $b \in B$. Since $B^o$ separates the points 
of $B$, we have 
$$\beta(a,b_1)b_2 = \beta(a,b_2)b_1 \iff
\forall \phi \in B^o, \ \phi(\beta(a,b_1)b_2) = \phi(\beta(a,b_2)b_1)$$
$$\iff   \forall \phi \in B^o, \ f_\beta(a)*\phi(b) = 
\phi*f_\beta(a)(b).$$
Hence $\beta(a,b_1)b_2 = \beta(a,b_2)b_1$ for any $a\in A$ and $b \in 
B$
if and only if $f_\beta(A) \subset Z(B^o)$. One also checks directly 
that
 $\beta(a_1,b)a_2 = \beta(a_2,b)a_1$ for any
$a\in A$ and $b\in B$ if and only if  
$f_\beta(a_1)\otimes a_2 = f_\beta(a_2) \otimes a_1$, for any $a\in A$.
Thus we have a bijection $\beta \mapsto f_\beta$
between $\mathcal{ZP}(A\otimes B)$ and
$\mathcal L(A,B^o)$. It is not hard to check that this bijection is a
group morphism.\hfill$\Box$

\medskip

\begin{rem}{\rm
It also is possible to give a Hopf algebra morphism
interpretation of general central $\bic$-pairings.}
\end{rem}

Let $A$ be a finite-dimensional Hopf algebra and let $D(A)$ be its 
Drinfeld
double. It is well known that $D(A)$ is a cocycle twist of 
$A \otimes (A^*)^{\rm cop}$ (\cite[Proposition 2.2, Remark
2.3]{doitake}). 
Combining the invariance of the second lazy cohomology group under 
twisting
(Corollary \ref{a-b}), the Schur-Yamazaki formula and
Lemma \ref{pairings} we get the following result.

\begin{coro}
Let $A$ be a finite-dimensional Hopf algebra and let $D(A)$ be its Drinfeld
double. We have a group isomorphism:
$$H^2_L(D(A)) \cong H^2_L(A) \times H^2_L(A^*) \times {\mathcal
L}(A).$$\hfill$\Box$
\end{coro}

\begin{example}{\rm Let $H_4$ be Sweedler's Hopf algebra. Its centre
is $k$. If $f\in{\mathcal L}(H_4)$, then $f\in {\rm Alg}(H_4,k)$ hence
$f=\varepsilon$ or $f=1^*-g^*$ with notation as in Example \ref{en}. 
Condition (\ref{ell}) is not verified for $f=1^*-g^*$ and $a=x$, so
${\mathcal L}(H_4)$ is trivial and by self-duality of $H_4$ it follows
that $H^2_L(D(H_4))\cong k\times k$.}
\end{example}

\section{Action on universal $r$-forms}\label{universal}

Let us suppose that $A$ is coquasitriangular with universal
$r$-form $r$. It is well known (\cite[Page 61]{mabook}) that if 
$\sigma$ is a left 2-cocycle
and if $\tau$ is the usual flip operator
then $(\sigma\tau)*r*\sigma^{-1}=r_\sigma$ is a universal $r$-form for 
the twisted Hopf algebra $_{\sigma}\!A_{\sigma^{-1}}$. In particular, 
if $\sigma$ is a lazy cocycle then $r_\sigma$ is again a universal
$r$-form for $A$ and this defines a left action of $\lazy$ on the set
$\U$ of universal $r$-forms of $A$. 

\begin{lemma}Let $A$ be a coquasitriangular Hopf algebra. Then
\begin{enumerate} 
\item The action of $\lazy$ on $\U$ factors through an action of 
$\coh2$.
\item The right action $\gets$ of $\aut$ on $\U$ induces an action of
${\rm CoOut}(A)$ on $\U$. 
\item The actions of $H^2_L(A)$ and of ${\rm CoOut}(A)$ on $\U$ combine
to a right action of ${\rm CoOut}(A) \ltimes H_L^2(A)$ given by 
$r\arp(\overline{\alpha},\overline{\sigma})=\sigma^{-1}\tau*(r\gets\overline{\alpha})*\sigma$.
\item The kernel of this action contains ${\rm CoOut}(A)^-$.
\end{enumerate}
\end{lemma} 
\pf {\it 1.} The 
group $B^2_{L}(A)$ acts trivially on $\U$. Indeed, let $\gamma\in
\Reg_L^1(A)$, 
let $r\in \U$ and let $a,\,b\in A$. Then
$$
\begin{array}{rl}
r_{\partial(\gamma)}(a,b)&=\partial(\gamma)(\b1,\a1) 
r(\a2,\b2)(\partial(\gamma))^{-1}(\a3,\b3)\\
&=\gamma(\a1)\gamma(\b1)\gamma^{-1}(\b2\a2)r(\a3,\b3)\gamma(\a4\b4)
\gamma^{-1}(\a5)\gamma^{-1}(\b5)\\  
&=\gamma(\a1)\gamma(\b1)r(\a2,\b2)\gamma^{-1}(\a3)
\gamma^{-1}(\k3)=r(a,b).\\
\end{array}
$$  
\noindent{\it 2.} It is well known that the right action $\gets$ of 
$\aut$ on
$\Reg^2(A)$ stabilizes $\U$ if $A$ is
coquasitriangular. If $\ad(\gamma)\in\CoInn(A)$ we may assume that
$\gamma$ is an algebra morphism $A\to k$. Then, for every
$a,b\in A$:
$$
\begin{array}{rl}
r(\ad(\gamma)(a),\ad(\gamma)(b))&=\gamma^{-1}(\a1)\gamma^{-1}(\b1)
r(\a2,\b2)\gamma(\a3\b3)\\
&=\gamma^{-1}(\b1)\gamma^{-1}(\a1)r(\a3,\b3)\gamma(\b2\a2)\\
&=r(a,b).
\end{array}
$$
\noindent{\it 3.} It is not hard to check that 
$r\arp(1,\overline{\sigma})(\overline{\alpha},1)=r\arp(\overline{\alpha},\overline{\sigma\gets\alpha})$.

\noindent{\it 4.} The group ${\rm CoOut}(A)^-$ is represented by pairs
$(\ad(\mu), \partial(\mu)^{-1})$. For every $a,b\in A$ we have:
$$
\begin{array}{l}
(r\arp(\ad(\mu), \partial(\mu)^{-1}))(a,b)\\
\phantom{r(a,b)}=\partial(\mu^{-1})^{-1}(\b1,\a1)r(\ad(\mu)(\a2),\ad(\mu)(\b2))
\partial(\mu)^{-1}(\a3,\b3)\\
\phantom{r(a,b)}=\mu^{-1}(\b1\a1)\mu(\a2)\mu(\b2)\mu^{-1}(\a3)
\mu^{-1}(\b3)r(\a4,\b4)\mu(\b5\a5)\\
\phantom{r(a,b)}=r(\a1,\b1)\mu^{-1}(\a2\b2)\mu(\a3\b3)=r(a,b).
\end{array}
$$\hfill$\Box$
% by $\alpha\lu
% \gamma=\gamma\circ(\alpha^{-1})^{\otimes q}$. This action stabilizes 
%$\Reg^q_L(A)$ because if
% $\gamma\in \Reg_L^q(A)$ then
%$$
%\begin{array}{l}
%\gamma(\alpha^{-1}(\a1^1)\otimes\cdots\otimes\alpha^{-1}(\a1^q))
%\a2^1\cdots \a2^q\\
%=\alpha\left(\gamma(\alpha^{-1}(a^1)_{(1)}\otimes\cdots\otimes\alpha^{-1}(a^q))_{(1)} 
%\alpha^{-1}(a^1)_{(2)}\cdots \alpha^{-1}(a^q)_{(2)}\right)\\
%=\alpha\left(\alpha^{-1}(a^1)_{(1)}\cdots \alpha^{-1}(a^q)_{(1)}
%\gamma(\alpha^{-1}(a^1)_{(2)}\otimes\cdots\otimes\alpha^{-1}(a^q)_{(2)})
%\right)\\ 
%=\a1^1\cdots \a1^q 
%\gamma(\alpha^{-1}(\a2)\otimes\cdots\otimes\alpha^{-1}(\a2^q)).
%\end{array}
%$$
%This action stabilizes the sets of coboundaries and the sets of
%cocycles. In particular, it stabilizes $\lazy$ and $B^2_L(A)$,
%inducing an action of $\aut$ on $\coh2$. If $A$ is dual
%quasitriangular, then the action of $\aut$ on $\reg2$ stabilizes
%$\U$. 
%If $r\in\U$ we shall denote by $r^\alpha$ the composition
%$r\circ(\alpha^{-1}\otimes\alpha^{-1})$.  

\medskip

Let $\cal G$ denote the semi-direct product of $\coh2$ and
$\aut$. Then 
$\cal G$ acts again on $\U$ and the $\cal G$-orbits on $\U$ coincide
with the ${\rm CoOut}(A) \ltimes H_L^2(A)$-orbits on $\U$. 
%
%$$
%\begin{array}{rl}
%(\sigma,\alpha)\lu r&=((\sigma,\,1)(1,\,\alpha))\lu r\\
%&=(\sigma,\,1)\lu r^\alpha\\
%&=(\sigma\circ\tau)*r^\alpha*\sigma^{-1}
%\end{array}
%$$
%and
%$$
%\begin{array}{rl}
%((1,\,\alpha)(\sigma,\,1))\lu r&=((\sigma\circ\tau)^\alpha*r^\alpha*(\sigma^{-1})^{\alpha}\\
%&=(\alpha\lu\sigma,\,\alpha)\lu r.
%\end{array}
%$$
%

\bigskip

Let $(A,r)$ be a coquasitriangular Hopf algebra. 
We recall that the Brauer group 
$BC(k,\,A,\,r)$ is the
Brauer group of the braided monoidal category of right $A$-comodules,
%where the algebras are the $A^{op}$-right
%comodule algebras and 
where the braiding is given on $M\otimes N$ by 
$\psi_{MN}(y\otimes n)=\n1\otimes y_1\,r(\n2\otimes y_2)$ if
the comodule
structures on $M$ and $N$ are given by $\rho_M(y)=y_1\otimes y_2$ and
$\rho_N(n)=\n1\otimes \n2$. 

In this new language, \cite[Proposition 3.1]{gio} becomes:

\begin{theorem}The Brauer group $BC(k,\,A,\,r)$ is constant on the
$H^2_L(A)$-orbits of $\U$.\hfill$\Box$   
\end{theorem}

\begin{corollary}The Brauer group $BC(k,\,A,\,r)$ is constant on the
${\cal G}$-orbits of $\U$. The $\cal G$-orbits on $\U$ coincide with 
the 
${\rm CoOut}(A)\ltimes H^2_L(A)$-orbits on $\U$. \hfill$\Box$
\end{corollary}

It is an interesting problem to understand the orbits under the $\cal 
G$-action.
Let us recall that if $(A,\,r)$ is cotriangular then $(A,\,s)$ will
be cotriangular for every $s\in r\arp {\cal G}$. 
One may wonder in which cases
the action of $\cal G$ on cotriangular structures is 
transitive. For instance, if $A=k[{\mathbb Z}_2]$ then there are 2
distinct universal $r$-forms on $A$: 
$\varepsilon\otimes\varepsilon$ and $\frac{1}{2}(1^*\otimes
1^*+g^*\otimes 1^*+g^*\otimes 1^*-g^*\otimes g^*)$ and they do not lie
in the same orbit because every cocycle for $k[{\mathbb Z}_2]$ is a 
symmetric form. On the other hand we have:

%
%\begin{example}In \cite{GioJuan1} it is shown for $H=H_\nu$ 
%that the orbits of $\U$ under the $H^2_L(H)$-action are parametrized
%by the equivalence classes of bicharacters
%for the cyclic group ${\mathbb Z}_{2\nu}$ 
%at least when $\nu\neq 4$, in particular the information is contained
%the restricion of the $r$-forms to the group of grouplike
%elements. This is also the parametrisation of the orbits under the
%$\cal G$-action. 
%The element $\sigma_{-1,0}$, which was not
%considered in \cite{GioJuan1} acts trivially on $\U$, so the above
%action is not faithful in general.
%\end{example}

\begin{example}{\rm Let $E(n)$ and $k$ be as in Example \ref{en}. 
By \cite{PVO1} the set of
quasitriangular structures for $E(n)$ (hence, dually, the set $\U$) is
parametrized by  $n\times n$ matrices with
entries in $k$ and the orbits under the action of
  ${\rm Aut}_{\rm Hopf}(E(n))$ on $\U$ correspond to congruence classes 
of
matrices. By 
\cite[Theorem 2.9]{GioJuan3} the orbits of the $\cal G$-action on $\U$
are parametrized 
by congruence classes of skew-symmetric $n\times n$ matrices with
coefficients in $k$, i.e., by the matrices of the form
$J_l=\left(
\begin{array}{ccc}
0_l&I_l&0\\
-I_l&0_l&0\\
0&0&0_{n-2l}
\end{array}\right)
$ for $0\le l\le\left[\frac{n}{2}\right]$.
All triangular $r$-forms lie in the orbit represented by $J_0=0$, 
and they correspond to all 
symmetric matrices. } 
\end{example} 

\section{A monoidal category of projective representations
for Hopf algebras}\label{projective} 

The Schur multiplier is the traditional
companion of the theory of projective representations of groups.
For Hopf algebras, there is still a naive
notion of a projective representation, that
is a representation of the  Galois object
corresponding to a 2-cocycle. However
since the convolution product of 2-cocycles
is no longer a cocycle in general, there
is no nice monoidal structure on the category
of such general projective representations.
However, when we restrict ourselves
to lazy 2-cocycles,
we are able to construct a monoidal
$G$-category (with $G=H_L^2(A)$) of projective representations of
$A$. Such categorical structures were considered by Turaev
in the setting of homotopy quantum field theory \cite{[Tu]}.
The base ring $k$ is a field in this section.

\smallskip
Let us first recall some notions introduced
in \cite{[Tu]}. We assume familiarity with monoidal categories.
First a  category
$\mathcal C$ is said to be \textbf{$k$-additive} 
if all the Hom's in $\mathcal C$ are $k$-modules 
and the composition of morphisms is bilinear over $k$. 

A $k$-additive category is said to be \textbf{$\Lambda$-split} 
for some set $\Lambda$ if
there exists a family $(\mathcal C_\lambda)_{\lambda \in \Lambda}$
of $k$-additive full subcategories of $\mathcal C$ such that
each object belongs to $\mathcal C_\lambda$
for some $\lambda \in \Lambda$, and such that 
for any pair of objects $V \in \mathcal C_\lambda$ and 
$W \in \mathcal C_\nu$ with $\lambda \not = \nu$, then
Hom$_\mathcal C(V,W) = \{0\}$
(in particular any non-zero object belongs to $\mathcal C_\lambda$
for a unique $\lambda \in \Lambda$).
The family of full subcategories 
$(\mathcal C_\lambda)_{\lambda \in \Lambda}$
is said to be a \textbf{$\Lambda$-splitting} of $\mathcal C$.

A \textbf{$k$-additive monoidal category} is a monoidal category
which is $k$-additive as a category, and such that the tensor
product is bilinear over $k$.

Let $G$ be a group. A \textbf{monoidal $G$-category over $k$}
is a $k$-additive monoidal category with left duality 
$\mathcal C$ which is $G$-split and such that

\noindent
(i) If $V \in \mathcal C_\lambda$ and
$W \in \mathcal C_\nu$, then $V \otimes W \in C_{\lambda \nu}$,

\noindent
(ii) If $V \in \mathcal C_\lambda$, then 
$V^* \in \mathcal C_{\lambda^{-1}}$.  

Note that $\textbf{1}$, the monoidal unit of $\mathcal C$,
 necessarily belongs to $\mathcal C_1$,
which is itself a $k$-additive monoidal category with left duality. 

\smallskip

We are going to construct, for a Hopf algebra $A$, a monoidal
$H^2_L(A)$-category over $k$, constisting of projective
representations of $A$.

\begin{defi}
Let $A$ be a Hopf algebra. A (finite-dimensional)
\textbf{projective representation of $A$} consists of a triplet
$(\sigma, V , \pi_V)$ where $\sigma \in Z^2_L(A)$ is a lazy
2-cocycle, $V$ is a finite dimensional vector space and
$\pi_V : A(\sigma) \longrightarrow {\rm End}(V)$ is an algebra
morphism.

Let $X= (\sigma, V , \pi_V)$ and $Y=(\omega,W,\pi_W)$ be some
projective representations of $A$. A \textbf{basic morphism} 
$f : X \longrightarrow Y$ is a linear map $f : V \rightarrow W$
such that there exists $\mu \in {\rm Reg}^1_L(A)$ satisfying
$$f \circ \pi_V(a) = \mu(a_1) \pi_W(a_2) \circ f,
\quad \forall a \in A.$$
A \textbf{morphism} $X\longrightarrow Y$ is a linear combination
of basic morphisms $X \longrightarrow Y$.
\end{defi}

It is easy to check that the composition
of two basic morphisms is again a basic morphism.
Therefore, extending the composition by bilinearity
to all morphisms, we obtain a 
\textbf{category of projective representations of $A$}, 
which we denote by $\mathbb P(A)$.
Clearly $\mathbb P(A)$ is a $k$-additive category.
We are going to show that $\mathbb P(A)$ has a structure
of monoidal $H^2_L(A)$-category over $k$.
The basic tool for proving the $H^2_L(A)$-splitting and for proving
the existence of left duals is the generalized antipode of the
Hopf-Galois system associated with a cocycle (\cite{[BiHo]}). 

\begin{lemm}\label{anti}
Let $\sigma$ be a left 2-cocycle on a Hopf algebra $A$.
Then the linear map 
$\phi_\sigma : \  _{\sigma} \! A \longrightarrow A_{{\sigma}^{-1}}$
defined by $\phi_\sigma(a) = \sigma(a_{1}, S(a_{2})) S(a_{3})$,
is an algebra anti-morphism.
Futhermore we have, for $a \in A$,
$$\phi_\sigma(a_1)._{\sigma^{-1}} a_2 = \varepsilon(a)1 =
a_1 ._{\sigma^{-1}} \phi_\sigma(a_2).$$ 
\end{lemm}
\pf Let $a,b \in A$. Since $\sigma$ is a left
2-cocycle, we get
$$
\begin{array}{l}
\sigma(a_1,b_1)\sigma(a_2b_2,S(a_3b_3))\\
\phantom{\sigma(a_1,b_1)}=\sigma(b_1,S(a_3b_4)) 
\sigma(a_1,b_2S(a_2b_3)) \\
\phantom{\sigma(a_1,b_1)}= \sigma(b_1,S(a_3b_2))\sigma(a_1,S(a_2)) \\
 \phantom{\sigma(a_1,b_1)}= \sigma^{-1}(S(b_5),S(a_4)) 
\sigma(b_1,S(b_4))
\sigma(b_2S(b_3), S(a_3)) \sigma(a_1,S(a_2)) \\
 \phantom{\sigma(a_1,b_1)} = \sigma(b_1,S(b_2))\sigma(a_1,S(a_2))
\sigma^{-1}(S(b_3),S(a_3)).
\end{array}
$$
Using this computation, it is immediate to see
that $\phi$ is an algebra anti-morphism. 
The last statement is proved in the proof of \cite[Proposition 
2.1]{[BiHo]} (In fact since $\phi_\sigma$ is the antipode
of the Hopf-Galois system associated with $\sigma$, the 
anti-multiplicativity of $\phi_\sigma$ also follows 
from \cite[Corollary 1.10]{[BiHo]}).\hfill 
$\Box$

\begin{lemm}
Let $X= (\sigma, V , \pi_V)$ and $Y=(\omega,W,\pi_W)$ be some
projective representations of $A$ and let 
$f : X \longrightarrow Y$ be a non-zero basic morphism.
Then $\overline{\sigma} = \overline{\omega}$ in $H^2_L(A)$.
\end{lemm}
\pf Let $\mu \in {\rm Reg}^1_L(A)$ satisfying
$f \circ \pi_V(a) = \mu(a_1) \pi_W(a_2) \circ f$,
$\forall a \in A$. Let $a,b \in A$. We have, using Lemma \ref{anti}
and the laziness of $\sigma$ and $\mu$,
$$
\begin{array}{rl}
\sigma(a_1,b_1)\mu(a_2b_2) f& 
=\sigma(a_1,b_1)\mu(a_2b_2) \varepsilon(a_3b_3)f \\
& =\sigma(a_1,b_1)\mu(a_2b_2) \pi_W(\phi_{\omega^{-1}}(a_3b_3)) \circ 
\pi_W(a_4b_4) \circ f \\ 
&= \pi_W(\phi_{\omega^{-1}}(a_1b_1)) \sigma(a_2,b_2) 
\mu(a_3b_3) \pi_W(a_4b_4) \circ f \\
&=\pi_W(\phi_{\omega^{-1}}(a_1b_1)) \sigma(a_2,b_2) f \circ
\pi_V(a_3b_3)\\
&= \pi_W(\phi_{\omega^{-1}}(a_1b_1)) \circ f \circ \pi_V(a_2) \circ 
\pi_V(b_2)\\
 &= \pi_W(\phi_{\omega^{-1}}(a_1b_1)) \mu(a_2) \mu(b_2) 
 \pi_W(a_3) \circ \pi_W(b_3) \circ f \\
 & =\pi_W(\phi_{\omega^{-1}}(a_1b_1)) \mu(a_2) \mu(b_2)
\omega(a_3,b_3) \pi_W(a_4b_4) \circ f  \\
&= \mu(a_1)\mu(b_1) \omega(a_2,b_2) 
\pi_W(\phi_{\omega^{-1}}(a_3b_3)) \pi_W(a_4b_4) \circ f\\
&=\mu(a_1)\mu(b_1) \omega(a_2,b_2) f.
\end{array}
$$
Hence, since $f$ is non-zero, we have $\sigma = \partial(\mu) \omega$.
\hfill$\Box$ 

\medskip

Let $x \in H^2_L(A)$. We define $\mathbb P_x(A)$ to be the full
subcategory of $\mathbb P(A)$ consisting of projective representations
$(\sigma,V,\pi_V)$ with $\overline{\sigma}=x$. It is clear from the 
above
lemma that the $k$-additive category $\mathbb P(A)$ is 
$H^2_L(A)$-split,
with $(\mathbb P_x(A))_{x \in H^2_L(A)}$ as $H^2_L(A)$-splitting.
 
\medskip

Let us now endow $\mathbb P(A)$ with a monoidal structure.
Let $X= (\sigma, V , \pi_V)$ and $Y=(\omega,W,\pi_W)$ be some 
projective representations of $A$. Since
$\Delta : A(\sigma * \omega) \rightarrow A(\sigma) \otimes A(\omega)$ 
is an algebra morphism, it follows that
$\pi_{V \otimes W} = (\pi_V \otimes \pi_W) \circ \Delta :
A(\sigma * \omega) \longrightarrow {\rm End}(V \otimes W)$,
is an algebra morphism. Thus we have defined a new projective 
representation of $A$:
$$X \otimes Y := (\sigma * \omega, V \otimes W, \pi_{V \otimes W}).$$
It is easy to check that the tensor product of two basic morphism
is still a basic morphism, and hence the tensor product of two
morphism is still a morphism. In this way $(\mathbb P(A), \otimes, 
\textbf{1})$,
endowed with the obvious associativity constraints and with
$\textbf{1}=(\varepsilon \otimes \varepsilon, k, \varepsilon)$ 
as monoidal unit, is a $k$-additive monoidal category. 
Let us now check that every object has a left dual.

\begin{lemm}
Let $X= (\sigma, V , \pi_V)$ be a projective representation of $A$.
Let $\pi_{V^*} : A \longrightarrow {\rm End}(V^*)$ be defined by
$\pi_{V^*}(a) = {^t \! \pi_V}(\phi_{\sigma^{-1}}(a))$.
Then 
$$X^* := (\sigma^{-1}, V^*, \pi_{V^*})$$
is a projective representation of $A$, and is a left dual for $X$. 
\end{lemm}
\pf It follows from Lemma \ref{anti}, that $\pi_{V^*}$
is an algebra morphism and hence $X^*$ is a projective representation.
Let $e : V^* \otimes V \rightarrow k$ be the evaluation map.
Let $f \in V^*$ and $v \in V$. Then, by Lemma \ref{anti}, we have 
for $a \in A$,
\begin{align*}
e  \circ \pi_{V^*\otimes V}(a)(f \otimes v) &=
e ((f \circ \pi_V(\phi_{\sigma^{-1}}(a_1)) \otimes \pi_V(a_2)(v))\\
& = f \circ  \pi_V(\phi_{\sigma^{-1}}(a_1)) \circ \pi_V(a_2) (v)\\
&= f \circ \pi_V(\phi_{\sigma^{-1}}(a_1)._{\sigma} a_2)(v) \\
& = \varepsilon(a) f(v) = \varepsilon(a) e(f \otimes v).
\end{align*}
This shows that $e : X^* \otimes X \rightarrow \textbf{1}$
is a morphism of projective representations.
One shows similarly
that the coevaluation map $\delta : k \rightarrow V \otimes V^*$
is a morphism $\textbf{1} \longrightarrow X \otimes X^*$.
Thus the triplet $(X^*,e,\delta)$ is a left dual for $X$
in $\mathbb P(A)$.

\phantom{.}\hfill$\Box$   

\medskip 
Let $X$ and $Y$ be some projective representations of $A$.
It is clear that if $X \in \mathbb P_x(A)$ and $Y \in \mathbb P_y(A)$,
then $X \otimes Y\in \mathbb P_{xy}(A)$ and $X^* \in \mathbb 
P_{x^{-1}}(A)$.
Summarizing the results of the section, we have proved:

\begin{theo}
Let $A$ be a Hopf algebra. Then $\mathbb P(A)$ is a monoidal
$H^2_L(A)$-category over $k$.\hfill$\Box$
\end{theo}

\section{Examples: monomial Hopf algebras}\label{monomial}

In this section $k$ is a field containing all primitive 
roots of unity.
We compute the second lazy cohomology group
for the monomial Hopf algebras. 

Recall (\cite{[CHYZ]}) that a \textbf{group datum} (over $k$) is a 
a quadruplet $\mathbb G = (G,g,\chi,\mu)$ consisting of
a finite group $G$, a central element $g \in G$,
a character $\chi : G \longrightarrow k^\cdot$ 
with $\chi(g) \not = 1$ and an element $\mu \in k$ such that
$\mu =0$ if $o(g) =o(\chi(g))$, and that if $\mu \not =0$,
then $\chi^{o(\chi(g))}=1$.

Let $\mathbb G = (G,g,\chi,\mu)$ be a group datum. 
A Hopf algebra $A(\mathbb G)$ is associated with $\mathbb G$ 
in \cite{[CHYZ]}. We will slightly change the conventions
of \cite{[CHYZ]} for the formula defining the coproduct
but this will not change the whole set of isomorphism classes.
As an algebra $A(\mathbb G)$ is the
quotient of the free product algebra
$k[x] * k[G]$ by the 
two-sided ideal generated by the relations
$$xh = \chi(h) hx, \ \forall h \in G , \quad \quad 
x^d = \mu(1-g^d), \ {\rm where} \ d = o(\chi(g)).$$  
The Hopf algebra structure of $A(\mathbb G)$ is defined
by:
$$\Delta(x) = 1  \otimes x + x \otimes g, \quad \varepsilon(x) =0, 
\quad
S(x) = -xg^{-1},$$
$$\Delta(h)= h \otimes h \ , \quad \varepsilon(h) =1 \ , \quad
S(h)=h^{-1} \ , \ \forall h \in G.$$ 
Using the diamond lemma \cite{[Be]}, it is not difficult
to see that the set 
$\{hx^i, 0 \leq i \leq d-1, h \in G \}$ is a linear basis
of $A(\mathbb G)$, and hence $\dim_k(A(\mathbb G))= |G|d$.

The Hopf algebras $A(\mathbb G)$ have been
shown in \cite{[CHYZ]} to be exactly the monomial non semisimple
Hopf algebras: see \cite{[CHYZ]} for the precise concept
of a monomial Hopf algebra. 

\medskip

We need the notion of type I group datum introduced in \cite{[BiG]}.
A \textbf{type I group datum} $\mathbb G = (G,g,\chi,\mu)$
is a group datum with $\mu =0$, $d= o(\chi(g)) = o(g)$
and $\chi^d =1$. 
In this case we simply write $\mathbb G = (G,g,\chi)$.
We will not need the other types of group data.
We can state now the main result of the section.

\begin{theo}\label{explicit}
Let $\mathbb G = (G,g, \chi, \mu)$ be a group datum. Then we have the 
following
description for $H^2_L(A(\mathbb G))$.
\begin{itemize}
\item If $\mathbb G$ is of type I, then
$H^2_L(A(\mathbb G)) \cong
H^2(G/\langle g \rangle,k^\cdot) \times k$.
\item If $\mathbb G$ is not of type I, then 
$H^2_L(A(\mathbb G)) \cong H^2(G/\langle g \rangle,k^\cdot)$.
\end{itemize}
\end{theo}

The proof will use the description of the biGalois objects given
in \cite{[BiG]}. Before going into the heart of the proof,
we need some cohomological preliminaries.
Let us first recall the definition, introduced in \cite{[BiG]}, 
of the modified second cohomology
group of a group.

Let $G$ be a group and let $g \in G$ be a central element. We put
$$Z^2_g(G,k^\cdot) = \{\sigma \in Z^2(G, k^\cdot), \ 
\sigma(g,h) = \sigma(h,g), \forall h \in G \} \quad {\rm and}$$  
$$B^2_g(G , k^\cdot) = \{ \partial(\mu), \
\mu : G \rightarrow k^\cdot, \ \mu(g)= 1 = \mu(1)\}.$$
For $g_1,g_2 \in Z(G)$, it is clear that
$B_{g_2}^2(G, k^\cdot)$ is a subgroup of $Z^2_{g_1}(G, k^\cdot)$ and
we define
$$H_{g_1,g_2}^2(G, k^\cdot) = Z^2_{g_1}(G,k^\cdot) / 
B_{g_2}^2(G,k^\cdot).$$
We have $ H_{1,1}^2(G, k^\cdot) =H^2(G, k^\cdot)$. 
Let us introduce another group now.
We define $Z^2_{L,g}(G,k^\cdot)$ to be the subset of $Z_g^2(G,k^\cdot)$
consisting of elements $\sigma$ such that there exists 
$\mu : G \longrightarrow k^\cdot$ satisfying $\mu(g) =1 = \mu(1)$
and $\sigma(g,h) = \mu(h) \mu(gh)^{-1}$, $\forall h \in G$.
It is clear that $Z^2_{L,g}(G,k^\cdot)$ is a subgroup of
$Z^2_g(G,k^\cdot)$ containing $B^2_g(G, k^\cdot)$. 
So we put 
$$L_g(G,k^\cdot) = Z^2_{L,g}(G,k^\cdot)/ B^2_g(G,k^\cdot).$$
This subgroup of $H^2_{g,g}(G, k^\cdot)$
will appear naturally in the study of the
bicleft biGalois objects. In fact we have:

\begin{lemm}\label{strange}
We have a group isomorphism
$L_g(G,k^\cdot) \cong H^2(G/ \langle g \rangle,k^\cdot)$.
\end{lemm}
\pf Let $\pi : G \longrightarrow G / \langle g \rangle$ be the
canonical surjection. Let $\sigma \in Z^2(G/\langle g\rangle, 
k^\cdot)$: it is
clear that $\sigma \circ (\pi \times \pi) \in Z^2_{L,g}(G,k^\cdot)$ and 
that
if $\sigma$ is a coboundary, then 
$\sigma \circ (\pi \times \pi)$ belongs to $B^2_g(G, k^\cdot)$. 
Therefore we get a group morphism
\begin{align*}
\theta : H^2(G/\langle g \rangle, k^\cdot) & \longrightarrow 
L_g(G,k^\cdot) \\
\overline{\sigma} & \longmapsto \overline{\sigma \circ (\pi \times 
\pi)}.
\end{align*}
Let us show that $\theta$ is an isomorphism.

\noindent
\textbf{Claim 1}. $\theta$ is injective.

\noindent
\textsl{Proof of Claim 1}. Let $\sigma \in Z^2(G/\langle g \rangle , 
k^\cdot)$
be such that there exists $\phi : G \rightarrow k^\cdot$
such that $\phi(g) = 1 = \phi(1)$ and $\sigma \circ (\pi \times \pi)
= \partial(\phi)$.
Computing $\sigma(\pi\times\pi)(g,h)=\sigma(\pi\times\pi)(1, h)=1$ 
we see that $\phi (hg) = \phi(h)$ for any $h \in G$.
Hence there exists $\phi' : G / \langle g \rangle \rightarrow
k^\cdot$ such that $\phi = \phi' \circ \pi$ and
then $\sigma \circ (\pi \times \pi) = \partial(\phi' \circ \pi)$.
We conclude that $\sigma$ is a coboundary and that
$\theta$ is injective. \hfill$\Box$

\noindent
\textbf{Claim 2}. Let $\mu : G\to k$
with $\mu(1) =1$
be such that $\partial(\mu) \in Z^2_{L,g}(G,k^\cdot)$.
Then $\overline{\partial(\mu)} \in {\rm Im}(\theta)$. 

\noindent
\textsl{Proof of Claim 2}. Let $\omega \in Z^2_{L,g}(G, k^\cdot)$. 
Then one can find $\omega' \in Z^2_{L,g}(G, k^\cdot)$
having the same class as $\omega$ in $L_g(G,k^\cdot)$ such that
$\omega'(g,h) = \omega'(h,g) = 1$. Hence we can assume
that $\mu(gh) = \mu(g) \mu(h)$, $\forall h \in G$.
With this assumption the restriction of $\mu$ to
$\langle g \rangle$ is a character.

Now let us fix $s: G/ \langle g \rangle \rightarrow G$ a section
of $\pi$ with $s(1)=1$. We have 
$G = \langle g \rangle s(G/\langle g \rangle)$ and we define a map
$\gamma : G \rightarrow k^\cdot$ by 
$\gamma(g^i s(X)) =\mu(s(X))$ for $i \in \mathbb Z$ and 
$X \in G/\langle g \rangle$. We have $\gamma(g)=1$.
Now consider the function
$f_s : G/\langle g \rangle \times G/\langle g \rangle
\rightarrow\langle g\rangle $ defined by 
$f_s(X,Y) = s(X)s(Y)s(XY)^{-1}$. Then $\mu \circ f_s$ is a 2-cocycle
(in fact the image of the character $\mu_{|\langle g \rangle}$
by the transgression map) and we put 
$\sigma = (\mu \circ f_s)^{-1}$. A straightforward computation shows
that $\partial(\mu) = (\partial(\gamma))*(\sigma\circ (\pi \times 
\pi))$,
and hence $\overline{\partial(\mu)} = \theta(\overline{\sigma})$.
\hfill$\Box$

\noindent
\textbf{Claim 3}. $\theta$ is surjective.

\noindent
\textsl{Proof of Claim 3}. Let $\sigma \in Z^2_{L,g}(G,k^\cdot)$.
As in the proof of Claim 2, we can assume
that $\sigma(g,h) = \sigma(h,g)=1$ for any $h  \in G$.
Then the restriction of $\sigma$ to $\langle g \rangle$
is trivial and
% considering the $2$-cocycle condition to $g,g^i,h$, we
%see that 
the group pairing $G \times \langle g \rangle
\rightarrow k^\cdot$, $(h,g^i) \mapsto \sigma(h,g^i) 
\sigma(g^i,h)^{-1}$,
is trivial. Hence we can use the exact sequence of 
Iwahori and Matsumoto (see \cite[Theorem 2.2.7]{[Kar]}):
there exists $\omega \in Z^2(G/\langle g\rangle,k^\cdot)$ and 
$\mu : G \rightarrow k^\cdot$
with $\mu(1)=1$ such that 
$\sigma = (\omega \circ (\pi \times \pi))* \partial(\mu)$.
Then  $\partial(\mu) \in Z^2_{L,g}(G,k^\cdot)$ and, by Claim 2,
$\overline{\sigma} \in {\rm Im}(\theta)$. \hfill$\Box$ 

The proof of the Lemma is now complete.\hfill $\Box$       

\medskip

Let $\mathbb G = (G,g, \chi, 0)$ be a group datum.
Let us recall now the description of the $A(\mathbb G)$-biGalois 
objects.
Let $\sigma \in Z^2(G,k^\cdot)$, let $u \in {\rm Aut}_g(G)$
(i.e. $u(g)=g$) and let $a \in k$. Assume that the triplet
$(\sigma,u,a)$ satisfies the following compatibility conditions:

\begin{equation}\label{compa1}\chi \circ u(h) = \sigma(g,h)^{-1}
\sigma(h,g) \chi(h),\quad  
\forall h \in G
\end{equation}

\begin{equation}\label{compa2} 
a=0 \ {\rm if} \ \mathbb G \ {\rm is} \  {\rm not}
\ {\rm of}\ {\rm type} \ {\rm I}
\end{equation}

Let the algebra $A_{\sigma, a}^u(\mathbb G)$ 
 be the algebra presented by generators 
$X,\,(T_h)_{h \in G}$ with defining relations,
$\forall h, h_1, h_2 \in G$:
$$T_{h_1}T_{h_2} = \sigma(h_1,h_2) T_{h_1h_2}, 
\quad T_1=1, \quad  XT_h = \chi(h) T_hX,
\quad X^d = a T_{g^d}.$$  
It is shown in \cite{[BiG]} that $A_{\sigma, a}^u(\mathbb G)$ 
is an $A(\mathbb G)$-biGalois object, with respective right and left
coactions $\rho : A_{\sigma,a}^u(\mathbb G)
\longrightarrow A_{\sigma,a}^u(\mathbb G) \otimes A(\mathbb G)$ 
and $\beta : A_{\sigma,a}^u(\mathbb G)\longrightarrow A(\mathbb G)
\otimes A_{\sigma,a}(\mathbb G)$ defined by 
$$\rho(X) = 1 \otimes x + X \otimes g,
\quad \rho(T_h) = T_h \otimes h, \quad \forall h \in G,$$
$$\beta(X) = 1 \otimes X + x \otimes T_g, \quad 
\beta(T_h) = u(h) \otimes T_h, \ \forall h \in G.$$
Every $A(\mathbb G)$-biGalois object is isomorphic to one of the form 
$A_{\sigma, a}^u(\mathbb G)$ for some triplet $(\sigma,u,a)$
satisfying (\ref{compa1}) and (\ref{compa2}). 

\begin{lemm}\label{triplet}
Let $(\sigma,u,a)$ be a triplet as above. Then the
$A(\mathbb G)$-biGalois object $A_{\sigma,a}^u(\mathbb G)$ is bicleft
if and only if $\sigma \in Z^2_{L,g}(G,k^\cdot)$ and $u = {\rm id}_G$.
\end{lemm}
\pf Let $f : A_{\sigma,a}^u(\mathbb G) \longrightarrow A(\mathbb G)$
be a bicolinear isomorphism. We can assume that $f(1)=1$.
By \cite[Proposition 2.3]{[BiG]} there is a right
$A(\mathbb G)$-colinear isomorphism $\Phi : A(\mathbb G) 
\longrightarrow A_{\sigma, a}^u(\mathbb G)$, $hx^i \longmapsto T_h 
X^i$.
Hence, since $f \circ \Phi$ is a right $A(\mathbb G)$-colinear 
automorphism
of $A(\mathbb G)$, there exists $\mu \in {\rm Reg}^1(A(\mathbb G))$ 
such that
$f \circ \Phi = \mu * {\rm id}$. Thus we have, for $h \in G$ and $0\leq 
i \leq d-1$
$$
f(T_hX^i) = \sum_{l=0}^i 
\begin{pmatrix}
i \\
l
\end{pmatrix}_{\!q}  
\mu(hx^{i-l}) hg^{i-l}x^l,
$$
where $q=\chi(g)$ and we have used the $q$-binomial coefficients.
In particular $f(T_h) = \mu(h)h$ for $h \in G$, and since $f$ is left 
colinear,
we find that $u = {\rm id_G}$. By condition (\ref{compa1}) 
$\sigma \in Z^2_g(G,k^\cdot)$. We also
have  $f(T_hX) = \mu(hx)hg + \mu(h)hx$ and using again the 
left colinearity of $f$, we find
that $\sigma(h,g) = \mu(h) \mu(hg)^{-1}$ and therefore $\mu(g)=1$, 
which means that $\sigma \in Z^2_{L,g}(G,k^\cdot)$. 
Conversely, assume that $\sigma \in Z^2_{L,g}(G,k^\cdot)$ and that
$u = {\rm id}_G$. 
As in the proof of claim 2 in the proof of 
Lemma \ref{strange} we can assume, without changing the class of 
$\sigma$
in $H^2_{g,g}(G,k^\cdot)$ and hence without changing the isomorphism
class of the $A(\mathbb G)$-bicomodule algebra $A_{\sigma,a}^u(\mathbb 
G)$
(\cite[Proposition 3.4]{[BiG]}), that 
$\sigma(g,h)=\sigma(h,g) = 1$, $\forall h \in G$.
Then define a linear isomorphism 
$f : A_{\sigma,a}^u(\mathbb G) \longrightarrow A(\mathbb G)$ by 
$f(T_hX^i) = hx^i$, for $h \in G$ and $0\leq i \leq d-1$.
One can check that $f$ is a bicolinear
isomorphism using $\sigma(g,h)=\sigma(h,g) = 1$, 
$\forall h \in G$ for the left colinearity.

\phantom{.}\hfill$\Box$

\medskip

\noindent\textbf{Proof of Theorem \ref{explicit}}. Let $\mathbb G =(G,g,\chi, 
\mu)$ be a group 
datum. Assume first that $\mathbb G$ is a type I group datum, i.e. that
$\mu=0$, that $o(\chi(g))=d=o(g)$ and that $\chi^d=1$. 
Then by Lemma \ref{triplet} we have a map
\begin{align*}
\Psi_0 : Z^2_{L,g}(G,k^\cdot) \times k & \longrightarrow {\rm 
Bicleft}(A(\mathbb G))\\
(\sigma,a) & \longmapsto A_{\sigma,a}^{\rm id}(\mathbb G).
\end{align*}
For $\sigma \in Z^2_{L,g}(G,k^\cdot)$, we have 
$\sigma(g,g) \ldots \sigma(g,g^{d-1}) =1$, so by
\cite[Proposition 3.5]{[BiG]} $\Psi_0$ is a group morphism.
Then $\Psi_0$ induces an injective group morphism
$\Psi : L_g(G,k^\cdot) \times k  \longrightarrow {\rm Bicleft}(A)$ by  
\cite[Proposition 3.4]{[BiG]}. 
Now let $Z$ be an $A(\mathbb G)$-biGalois object: by
\cite[Proposition 3.8]{[BiG]} there exists a triplet
$(\sigma,u,a)$ as above such that $[Z] = [A_{\sigma, a}^u(\mathbb G)]$.
If $Z$ is bicleft, we have $\sigma \in Z^2_{L,g}(G,k^\cdot)$ and 
$u = {\rm id}_G$ by Lemma \ref{triplet} and therefore $\Psi$ is 
surjective and
it is an isomorphism. 
Then Lemma \ref{strange} concludes the proof in the type I case.

Assume now that $\mathbb G = (G,g,\chi,0)$ is not of type I.
Then the proof, using the results in \cite[Section 3]{[BiG]},
 is essentially the same as the one of the type I case. 
This is left to the reader.
%Assume now that $\mathbb G = (G,g,\chi,0)$ is not of type I.
%Then the proof, using the results in \cite[Section 3]{[BiG]},
% is essentially the same as the one of the type I case,
%the only difference being that $A_{\sigma, a}^u(\mathbb G)$ is 
%$A(\mathbb G)$-biGalois if and only if $a=0$. This is left to the 
%reader.
 
Finally if $\mathbb G = (G,g,\chi,\mu)$ with $\mu \not = 0$, then by 
\cite[Corollary 3.18]{[BiG]}, there exists an 
$A(\mathbb G)$-$A(\mathbb G_{\rm red})$-biGalois object 
for $\mathbb G_{\rm red} = (G,g,\chi,0)$ so the statement follows from
the previous case and Corollary \ref{a-b}\hfill$\Box$

\begin{ex}{\rm Recall (\cite{[BiG]}) that a cyclic datum
is a datum 
$(d,n,N,\alpha,q)$ where $d,n,N>1$ are integers, $\alpha \in \mathbb 
N^\cdot$
and $q \in k^\cdot$ is a root of unity, satisfying:
$$d|n|N, \ \alpha | \frac{N}{n}, \ {\rm GCD}(\alpha,d)=1, \ o(q) = 
\frac{Nd}{\alpha n}.$$ 
To any cyclic datum $(d,n,N,\alpha,q)$
we associate a group datum
$$\mathbb C[d,n,N,\alpha,q] :=
(C_N=\langle z \ | \ z^N=1 \rangle, g = z^{\frac{N}{n}}, \chi_q,0)$$
where $\chi_q$ is the character defined by $\chi_q(z) = q$. 
The associated Hopf algebra 
$A(d,n,N,\alpha,q) = A(\mathbb C[d,n,N,\alpha,q])$ is then the 
algebra presented by generators $z,x$ submitted to the relations
$$x^d=0, \quad z^N=1, \quad xz = qzx.$$
The coproduct is defined by $\Delta(z) = z \otimes z$ and 
$\Delta(x) = 1 \otimes x + x \otimes z^{\frac{N}{n}}$.
Using \cite[Lemma 4.4]{[BiG]} and Theorem \ref{explicit} we find that
$$H^2_L(A(d,n,N,\alpha,q)) \cong \left\{
\begin{array}{ll}
k &\mbox{ if } d=n =N,\\
(k^\cdot/(k^\cdot)^{\frac{N}{n}}) \times k & \mbox{ if } 
d=n<N , \\
& {\rm GCD}(\frac{N}{n},n) =1
\mbox{ and }\alpha = \frac{N}{n},\\
k^\cdot /(k^\cdot)^{\frac{N}{n}}
&\mbox {\rm otherwise}.
\end{array}\right.
$$
As a particular case, for the Taft algebras $H_{N,q} = A(\mathbb 
C[N,N,N,1,q])$,
we get $H^2_L(H_{N,q}) \cong k$.}
\end{ex}

\begin{ex}\label{counterex}{\rm Let 
$$G=\langle a,b,g~|~a^2=1=b^2=g^4,\; ag=ga,\; bg=gb,\;
    ab=bag^2\rangle.$$
$G$ is a non-abelian group of order $16$ with $g\in Z(G)$ and
    $o(g)=4$. Let $\chi\colon G\to k^\cdot$ be the character defined
    by $\chi(a)=\chi(b)=1$, $\chi(g)=-1$. We consider the group datum
    ${\mathbb G}=(G, g,\chi)$. For any $\phi\in{\rm Alg}(A({\mathbb
    G}), k)$ we have $\ad(\phi)(x)=\pm x$. Let $\mu_0\colon G\to
    k^\cdot$ be defined by: $\mu_0(a^\alpha b^\beta g^\gamma)=
(\sqrt{-1})^\gamma$. Now let $\mu\in{\rm Reg^1}(A({\mathbb G}))$ be
    defined by: $\mu(x^ih)=\delta_{i, 0}\mu_0(h)$ for $h\in G$. Since
    $\mu_0(gh)=\mu_0(g)\mu_0(h)$ for every $h\in G$ 
it follows that $\ad(\mu)$ is a Hopf algebra automorphism of
    $A({\mathbb G})$ with $\ad(\mu)(x)=\sqrt{-1}x$. Therefore
    $\ad(\mu)$ is not coinner, so $\CoInt(A({\mathbb G}))\neq
\CoInn(A({\mathbb G}))$ which means that 
${\rm CoOut}^-(A({\mathbb G}))$ is non-trivial.}\end{ex} 

\section{Lazy cohomology for some cotriangular Hopf algebras}\label{special}

In this section the base field will be $\mathbb C$. We shall start
recalling notation and results in \cite{chev1}, \cite{geleti1},
\cite{geleti2} and \cite{geleti3}. For terminology we mainly refer to
these papers.

\bigskip

Let $\cal A$ be a finite-dimensional Hopf superalgebra, with a 
grouplike
element $g$ such that $gxg^{-1}=(-1)^{\deg(x)} x$ for every
homogeneous element in ${\cal A}$. Then we can apply bosonization (see
\cite[\S 9.4 ]{mabook}, \cite[\S 3.1]{chev1}) obtaining a 
finite-dimensional
Hopf algebra $A$. The Hopf algebra $A$ is 
equal to $\cal A$ as an algebra but with coproduct
given, for homogeneous elements, by: 
$\Delta_A(h)=g^{\deg(\h2)}\h1\otimes(-1)^{(\deg(\h2)(1+\deg(h)}\h2$ if $\Delta_{\cal
  A}(h)=\h1\otimes\h2$. The Hopf superalgebra ${\cal A}$ is triangular
with (necessarily even) $R$-matrix ${\cal R}={\cal R}_0+{\cal R}_1\in
{\cal A}_0\otimes{\cal A}_0+{\cal A}_1\otimes {\cal A}_1$ if and only 
if $A$ is 
triangular with
$R$-matrix $R=({\cal R}_0+(1\otimes g){\cal R}_1)R_g$, where 
$R_g=\frac{1}{2}(1\otimes
1+1\otimes g+g\otimes 1-g\otimes g)$.   

Let $\cal A$ be a finite-dimensional cocommutative Hopf superalgebra
over  $\mathbb C$. By
\cite[Theorem 3.3]{kostant} 
${\cal A}={\mathbb C}[G]\ltimes \wedge W$ where $W$ is the
(purely odd) 
space of primitive elements and $G$ is the group of grouplikes
acting on $W$, hence on $\wedge W$, by conjugation.

If $G$ contains an element such that $g^2=1$ and such that
$gxg^{-1}=(-1)^{\deg(x)}x$, the procedure above
described yields a triangular Hopf algebra $A$ with $R$-matrix 
$R_g$.
 
By \cite[Proposition 3.4.1]{chev1} if $r\in S^2(W)$, the symmetric
algebra of $W$ (viewed in $W\otimes W$), 
then ${\cal J}=e^{r/2}\in{\cal A}\otimes{\cal A}$ satisfies
\begin{equation}\label{twist}(\Delta\otimes\id_A)({\cal J}){\cal
J}_{12}=(\id_A\otimes\Delta)({\cal J}){\cal
J}_{23}\mbox{ and }(\varepsilon\otimes\id_A)({\cal 
J})=(\id_A\otimes\varepsilon)({\cal J})=1,
\end{equation}
i.e, it is
{\bf a Drinfeld twist}  for ${\cal
  A}$. If we write ${\cal J}={\cal J}_{0}+{\cal J}_{1}$ with ${\cal 
J}_i\in
{\cal A}_i\otimes{\cal A}_i$ for $i=0,1$, then $J={\cal 
J}_{0}-(g\otimes
1){\cal }J_{1}$ is a Drinfeld twist for $A$, i.e., $J$ is a right
$2$-cocycle for $A^*$. It is not hard to check
that if we take $r\in S^2(W)^G$, the invariants under the $G$-action, 
then $J$ commutes with $\Delta(a)$ for every $a\in A$, i.e, $J$ is a 
lazy
cocycle for $A^*$. We shall call a lazy cocycle for $A^*$ also a 
{\bf lazy twist} for $A$. The dual version of
cohomology of cocycles is {\bf gauge equivalence}: 
two Drinfeld twists $J$ and $F$ for a Hopf
algebra $H$  
are said to be {\bf gauge 
equivalent} if $F=\Delta(x)J(x^{-1}\otimes x^{-1})$ for some $x\in H$
(see \cite{mabook}, \cite{geleti1} for details). In particular, two
gauge equivalent Drinfeld twists $F$ and $J$ are
cohomologous in lazy cohomology for $H^*$ if and only if the
element $x$ can be chosen to be central in $H$.

\medskip

By \cite{chev1}, \cite{geleti2} the triangular Hopf algebra $A$, 
with $R$-matrix $R_g$ 
is the key model of finite-dimensional triangular Hopf algebras over 
$\mathbb
C$. Indeed, all other such Hopf algebras are obtained from a Hopf
algebra of this type twisting the coproduct. In other words, for every
finite-dimensional triangular Hopf algebra $H$ there exists a Drinfeld
twist $J$ such that $H^J\cong A$ for some $G$ and $W$, where $H^J$ has
the same underlying algebra as
$H$ and $\Delta_{H^J}(h)=J^{-1}\Delta_H(h)J$ for every $h\in H$. The
$R$-matrix of $H^J$ is $(\tau J)^{-1}R_gJ$.   
Therefore, 
$A^*$ is the key model of finite-dimensional cotriangular Hopf
algebras over $\mathbb C$ and we have a method for the construction
of special lazy $2$-cocycles. More precisely,
\begin{lemma}\label{assignment}Let $A=({\mathbb C}[G]\ltimes\wedge 
W)^*$ with $g\in A$ and
  with coproduct
  as before. Then the assignment $r\mapsto {\cal J}=e^{r}\mapsto J$ 
defines an injective group morphism
$\Sigma\colon S^2(W)^G\to H^2_L(A)$.
\end{lemma}
\pf If ${\cal J}=e^{m}$ and ${\cal J}'=e^{m'}$ in the corresponding
Hopf superalgebra, then
$J=1\otimes 1-(g\otimes 1)m+\cdots$ and $J'=1\otimes 1-(g\otimes
1)m'+\cdots$ so $J*J'=1\otimes1-(g\otimes 1)(m+m')+\cdots$ and the
assignment gives a group morphism $S^2(W)^G\to Z^2_L(A)$. 
Combined with the standard projection we have a group morphism
$\Sigma\colon S^2(W)^G\to H^2_L(A)$. If for some $r\in S^2(W)^G$ we had 
$\Sigma(r)=(z^{-1}\otimes z^{-1})\Delta(z)$ for some central $z\in
A^*$, then the $R$-matrix of $A^*$ obtained through $\Sigma(r)$ would 
be
$\tau(\Delta z^{-1})R_g\Delta(z)=R_g$ because $R_g$ is an $R$-matrix. 
The corresponding $R$-matrix in
the Hopf superalgebra ${\cal A}^*$ is $e^{2r}=R_g^2=1\otimes 1$. Hence
$r=0$ and $\Sigma$ is injective.\hfill$\Box$

\medskip

Example \ref{en} shows that for the
family  $E(n)$ this construction exhausts all
the lazy $2$-cohomology classes. This holds in a more general 
framework. 

%Let us recall the classification in \cite{geleti1}. Every
%finite-dimensional triangular Hopf  algebra is a Drinfeld 
%twist of a triangular Hopf algebra with minimal part of dimension
%$\le2$. Such Hopf algebras are called modified supergroup algebras
%because they can be obtained ....
%
Let $(A,R)$ be a finite-dimensional cotriangular Hopf algebra. Then
there exists a cocycle $\omega$, a group $G$ acting on a vector space
$W$, and a central element $g\in
G\subset A^*$ acting as $-1$ on $W$ and such that $g^2=1$ 
for which
$(_\omega A_{\omega^{-1}})^*\cong {\mathbb C}[G]\ltimes\wedge W$ and
the $r$-form corresponding to $R$ under twist is
% either $\varepsilon\otimes\varepsilon$ or
$R_g:=\frac{1}{2}(\varepsilon\otimes\varepsilon+\varepsilon\otimes
g+g\otimes\varepsilon-g\otimes g)$. In particular, if $g=\varepsilon$
then $W$ is trivial and $R_g=\varepsilon\otimes\varepsilon$.
The data $G$ and $W$ are
unique up to isomorphism and the twist is unique up to gauge
equivalence.
By abuse of language, we will also say that $G$ and $W$ are data
associated to $A$. 

\medskip

Let $A$ be a cotriangular Hopf algebra with associated data $G$ and
$W$. We shall denote by $\{w_1,\cdots, w_n\}$ a fixed basis of $W$ and 
we shall denote by $\rho\colon G\to {\rm GL}_{n}({\mathbb C})$ 
the group morphism given by: $g^{-1}w_ig=\sum_j\rho(g)_{ij}w_j$ for 
every $i$. 

\begin{theorem}\label{triangular}
Let $A$ be a finite-dimensional cotriangular Hopf algebra, with
associated data $G$ and  $W$. 
If the representation $\rho$ of $G$ on $W$ is faithful then
$H^2_L(A)\cong S^2(W)^G$. 
\end{theorem}
\pf By Corollary \ref{a-b} it is enough to prove the result
when $A^*={\mathbb C}[G]\ltimes\wedge W$ because the hypothesis on
$\rho$ still holds if we twist the coproduct of $A^*$.
We shall use the terminology of Drinfeld twists rather than the
terminology of $2$-cocycles. 
Let $F$ be a lazy twist
for $A^*$ and let $\pi$ be the Hopf algebra projection of $A^*$
onto ${\mathbb C}[G]$. Then $\overline{F}=
(\pi\otimes\pi)(F)$ is a lazy twist for 
${\mathbb C}[G]$ and, since ${\mathbb C}[G]$ is also a sub Hopf algebra 
of
$A^*$, 
it is a twist for $A^*$. Since $F$ commutes with
$g\otimes g$, it is even in the ${\mathbb Z}_2$-gradation
induced by the 
action of $g$. Then, if ${\rm Rad}(A^*)$ denotes the Jacobson radical
  of $A^*$, i.e., the ideal generated by the $w_j$'s,
 $F=\overline{F}+$ terms in $({\rm Rad} (A^*))^2\otimes
A^*+A^*\otimes({\rm Rad}(A^*))^2+{\rm Rad}(A^*)\otimes{\rm Rad}(A^*)$.  
The elements of the form $h w_{i_1}\cdots
w_{i_m}$ with $1\le i_1<\cdots<i_m\le n$ and $h\in G$ form a basis for
$A^*$. Looking at the the expression of
$F\Delta(w_i)=\Delta(w_i)F$ as a linear combination of the
corresponding basis of $A^*\otimes A^*$ we see that
also $\overline{F}$ commutes with $\Delta(w_i)$ for every $i$ 
(hence it is a lazy twist for $A^*$). 
This implies that if ${\overline F}=\sum_{s,h\in G}f_{sh}s\otimes h$ then
$$\sum_{s,h\in G}f_{sh}sg\otimes hw_i=\sum_{l,p\in G}\sum_{j}
f_{lp}lu\otimes p\rho(p)_{ij}w_j$$
and
$$\sum_{s,h\in G}f_{sh}sw_i\otimes h=\sum_{l,p\in G}\sum_j
f_{lp}l\rho(l)_{ij}w_j\otimes p.$$
The coefficient of $gs\otimes hw_j$ in the first equality is: 
$\delta_{ij}f_{sh}=f_{sh}\rho(h)_{ij}$ so if $h$ is
such that $f_{sh}\neq0$ for some $s\in G$, then $h=1$ because $\rho$
is injective. With a similar
computation from the second equality we get that
$\overline{F}=1\otimes 1$.  

%If $F$ is a lazy coboundary, then $\overline F$ is
%a lazy coboundary for ${\mathbb C}[G]$. Therefore the projection 
defines
Let us consider $K=((A^*)^F)$ and its corresponding
Hopf superalgebra ${\cal K}=(({\cal A}^*)^{\cal
  F})$. 
The minimal part of $\cal K$ satisfies the conditions of \cite[Lemma
  5.3.2]{chev1}, hence its $R$-matrix is of the form 
${\cal R}_{\cal K}=e^r$ for $r\in
S^2(W)$. Since the corresponding $R$-matrix  ${\cal R}_{\cal K}$
commutes with $\Delta(w)$ for every $w$ in $W$, $r\in S^2(W)^G$. By
the classification in \cite{geleti2} 
$\cal K$ is isomorphic, as a triangular Hopf superalgebra, to $({\cal
  A}^*)^{\cal J}$ where ${\cal J}=e^{r/2}$. It follows from the 
classification in
\cite{geleti1} that the correspondig twist $J$ of $A^*$ is
gauge equivalent to $F$, i.e., $F=\Delta(z)*J*(z^{-1}\otimes z^{-1})=
J*\Delta(z)*(z^{-1}\otimes z^{-1})$ for $z\in A^*$ with
$\Delta(z)*(z^{-1}\otimes z^{-1})$ 
a coboundary for $A$ commuting with $\Delta(A^*)$.  
Since $(\pi\otimes\pi)(J)=1\otimes1$, the projection of $z$ is
grouplike, so
$\pi(z)=h$ for some $h\in G$. If we replace $z$ by $zh^{-1}$ we see
that $F=J*\Delta(zh^{-1})*(hz^{-1}\otimes hz^{-1})$
with $\nu=zh^{-1}-1$ nilpotent. 
$\Delta(1+\nu)*((1+\nu)^{-1}\otimes (1+\nu)^{-1})$ centralizes 
$\Delta(A^*)$
if and only if $l\mapsto (1+\nu)^{-1}l(1+\nu)$ is a coalgebra map in
$A^*$. But then for $h\in G$,  $(1+\nu)^{-1}h(1+\nu)$ is grouplike if 
and
only if it coincides with $h$ and for $w\in W$,  $(1+\nu)^{-1}w(1+\nu)$ 
is $(g,1)$ skew-primitive if and
only if it coincides with $w$. Hence $1+\nu$ is central in $A^*$ so 
$\Delta(1+\nu)*((1+\nu)^{-1}\otimes (1+\nu)^{-1})\in B^2_L(A)$. 
So if $F$ is a lazy $2$-cocycle,
then $F$ is cohomologous to $J$ with ${\cal J}=e^m$ and $m\in
S^2(W)^G$. 
%Since we already know that every such a $J$ is a lazy
%$2$-cocycle for $A$, we have a map $\Sigma\colon S^2(W)^G\to
%H^2_L(A)$. 
It follows
%from the previous discussion 
that the map $\Sigma$ in Lemma \ref{assignment} is surjective, whence 
the proof. 
%If ${\cal J}=e^{m}$ and ${\cal J}'=e^{m'}$ then
%$J=1\otimes 1-(u\otimes 1)m+\cdots$ and $J'=1\otimes 1-(u\otimes
%1)m'+\cdots$ so $J*J'=1\otimes1-(u\otimes 1)(m+m')+\cdots$ so 
%$\Sigma$
%is a group morphism. If for some $r\in S^2(W)^G$ 
%$\Sigma(r)=(z^{-1}\otimes z^{-1})\Delta(z)$ for some central $z\in
%A^*$, then the $R$-matrix of $A^*$ obtained through $\Sigma(r)$ is
%$\tau(\Delta z^{-1})R_u\Delta(z)=R_u$. The corresponding $R$-matrix in
%the Hopf superalgebra ${\cal A}^*$ is $e^{2r}=R_u^2=1\otimes 1$. Hence
%$r=0$ and $\Sigma$ is a group isomorphism.
\hfill$\Box$ 
%  
%we have the exact seuqence 
%\begin{equation}
%1\longrightarrow S^2(W)^G\longrightarrow H^2_L(A,{\mathbb C})
%\longrightarrow H^2_L(({\mathbb C}G)^*,{\mathbb C}).
%\end{equation} 
%\noindent{\bf Claim 2}The map $p$ is trivial.
%
%\noindent{\em Proof of Claim 2} .

\medskip

Let us observe that in this case $H^2_L(A)$ coincides with
the Hochschild cohomology of $A^*$, which is computed in \cite[Lemma
3.2]{geleti2}. 

\begin{remark}{\rm Let us observe that a Hopf algebra 
$A$ isomorphic to ${\mathbb C}[G]\ltimes \wedge W$ with $g$ as before
and with faithful
$G$-action is in general not self dual, even for $G$ abelian. 
Indeed, the intersection of the centre
of $A$ with its grouplikes is trivial, while the intersection of the
grouplikes of its dual with the centre is ${\rm Alg}({\mathbb
C}[G/\langle g\rangle],k)$.}
\end{remark}

\begin{proposition}
Let $A=({\mathbb C}[G]\ltimes\wedge W)^*$. 
Suppose that there exists  $g\in G$ central
acting as $-1$ on $W$ and
such that $g^2=1$  and suppose that the $G$-action on $W$ is faithful. 
Then 
$\CoInn(A)=\CoInt(A)$.  
\end{proposition}
\pf An invertible element $x$ in $A^*={\mathbb C}[G]\ltimes\wedge W$ is 
in $\Reg^1_{aL}(A)$ if and
only if conjugation by $x$ in $A^*$ is a coalgebra morphism. If $x$ is
central, there is nothing to prove. If $x$ is
not central, then $x=\pi(x)+y=\pi(x)(1+\nu)$ for some $\nu,y\in {\rm
  Rad}(A^*)$, with $\pi(x)\in{\mathbb C}[G]$ invertible. 
Conjugation by $\pi(x)$ in ${\mathbb C}[G]$ is a coalgebra map. 
Let $x=x_0+x_1$ in the ${\mathbb Z}_2$-gradation induced by conjugation
by $g$. Then $x^{-1}gx=h\in G$ and $x_0=x-x_1$ is invertible because
all elements in the radical of $A^*$ are nilpotent. Hence,
$x_0h+x_1h=gx_0+gx_1=x_0g-x_1g$ so $x_0hg=x_0$ obtaining that
$h=g$ and $x_1=0$. In particular,  $x=x_0$ is even and conjugation by 
$x$
maps each $w_i$ to a linear combination of the $w_j$'s:
$x^{-1}w_ix=(1+\nu)^{-1}\pi(x)^{-1}w_i\pi(x)(1+\nu)=\sum_j t_{ij}w_j$ 
and
$$\pi(x)(1+\nu)\sum_j t_{ij}w_j=w_i\pi(x)(1+\nu).$$ 
Looking at linear combinations of elements of the basis of $A^*$ we
see that this implies that conjugation by $\pi(x)$ gives already a
coalgebra morphism. 
Then $w_i\pi(x)=\sum_j t_{ij}\pi(x)w_j$. Putting $\pi(x)=\sum_{h\in 
G}c_hh$,
we have
$$\sum_{h\in G}c_hw_ih=\sum_j\sum_{v\in G}
t_{ij}c_vvw_j=\sum_{j,l}\sum_{v\in G}t_{ij}c_v\rho(v^{-1})_{jl}w_lv.$$
As in the proof of Theorem \ref{triangular}, 
we have that
$\delta_{il}c_h=c_h\sum_j t_{ij}\rho(h^{-1})_{jl}$ for every $i,l=1,\cdots
n$ and every $h\in G$. Then if $c_h,c_v\neq0$, $\rho(h)=\rho(v)=T=(t_{ij})$ 
so
by faithfulness of $\rho$, $h=v$ and $\pi(x)=c_vv$ for some $v$. Up to
rescaling of $x$, conjugation by $\pi(x)$ gives an element of
$\CoInn(A)$, so we might as well assume that $x=1+\nu$. Again, using 
the
basis of $A^*$ we see that conjugation by $1+\nu$ can be a coalgebra
morphism if an only if $\nu$ is central, hence the 
statement.\hfill$\Box$  

%This is a consequence of the proof of Theorem \ref{triangular}.  
%\hfill$\Box$

\begin{theorem}\label{triangular-trivial}
Let $A$ be a complex, finite-dimensional cotriangular Hopf algebra, 
with
associated data $G$ and $W$ and with central element $g\in A^*$. 
If $G=G'\times \langle g\rangle$ and the representation of $G'$ on $W$ 
is trivial 
then
$H^2_L(A)\cong S^2(W)\times H^2_L(({\mathbb C}[G'])^*)$. 
\end{theorem}
\pf As before it is enough to show the statement for $A=({\mathbb
C}[G]\ltimes \wedge W)^*$. If the representation of $G'$ on $W$ is 
trivial then $A^*\cong
{\mathbb C}[G']\otimes ({\mathbb C}[g]\ltimes \wedge W)$, so
$A\cong ({\mathbb C}[G'])^*\otimes E(n)$ where $n=\dim W$. 
By Theorem \ref{Yamazaki} 
$$H^2_L(A)\cong H^2_L({\mathbb
C}[G']^*)\times H^2_L(E(n))\times{\cal ZP}(({\mathbb
C}[G'])^*, E(n)).$$
By Example \ref{en}, it is enough to show that ${\cal ZP}({\mathbb
C}[G']^*,E(n))$ is trivial. Let $B$ be any Hopf algebra and let
$\beta$ be a central pairing between $B$ and $E(n)$. 
We have, for $b\in B$ and $w\in W$ and for $a=0,1$:
$$\beta(b,g^aw)g^a+\beta(b,g^{a+1})g^aw=
\beta(b,g^a)g^aw+\beta(b,g^aw)g^{a+1}.$$ 
Therefore $\beta(b, g^aw)=0$ and $\beta(b,g^a)=\beta(b, g^{a+1})$ 
for every $b\in B$, avery $a=0,1$ and every $w\in W$.  
Since $\beta(b,1)=\varepsilon(b)$, we see that $\beta(b,l
)=\varepsilon(b)\varepsilon(l)$ for every $b\in B$ and
every $l\in{\mathbb C}[{\mathbb
Z}_2]\cdot W$. Since $\beta(b,cc')=\beta(\b1,c)\beta(\b2,c')$, 
the group ${\cal ZP}(B, E(n))$ is always trivial and
we have
the statement.\hfill$\Box$

\bigskip

%Since
%$\overline{F}$ commutes with $\Delta(g)$ for every $g\in G$, it
%follows that $g\mapsto x^{-1}gx$ defines a group morphism, which
%is necessarily the identity. 
%
%Therefore, the
%categories of representations of $G$ with braiding given either by
%$r_u$ or by ${\overline F}_{21}^{-1}{\overline F}r_u$ are
%equivalent. 
%Since $u$ is central, these categories are equivalent if
%and only if the categories with braiding $1\otimes 1$ or
% ${\overline F}_{21}^{-1}{\overline F}$ are
%equivalent.  
%
%Since the categories   
%
%
%
%
%If  $(S^2V)^G=0$, then \cite[Proposition
%  3.3]{geleti2} shows that the there exist only finitely many 
%$2$-cocycles
%  $\sigma$, up to coboundaries, such that $_\sigma
%  H_{\sigma^{-1}}\cong H$, so $H^2_L(A)$ is a finite group. 
%It
%would be interesting to see under which assumptions $H^2_L(A)$ is
%trivial in this case and, in general, under which assumptions
%$H^2_L(A)=(S^2V)^G$ for this special class of Hopf algebras over
%$\mathbb C$. 
%
%WHAT HAPPENS IN THE GENERAL CASE
%
%DUAL PICTURE
%

\section{Appendix: Laziness and general Hopf-Galois extensions}

In this Appendix we study the relations between lazy cocycles and
general
crossed systems. 
 
The actions in Remark \ref{crossed} extend to all crossed
systems and are not defined only on cocycles with values in $k$. 
We recall that a crossed system over a $k$-algebra $R$ is 
a pair $(\lu, \sigma)$ where $\lu\colon A\to {\rm End}(R)$
 is a measuring of $A$ on $R$, 
$\sigma$ is a convolution invertible linear map $A\otimes A\to R$ and
 they satisfy the relations:
$$
\begin{array}{l}
\sigma(a, 1)=\sigma(1, a)=\varepsilon(a)1\\
(\a1\lu\b1\lu x)\sigma(\a2,\b2)=\sigma(\a1,\b1)(\a2\b2\lu x)\\
\sigma(\a1,\b1)\sigma(\a2\b2, 
c)=(\a1\lu\sigma(\b1,\c1))\sigma(\a2,\b2\c2)
\end{array}
$$
for every $x\in R$ and for every $a,\,b,\,c\in A$. 

\begin{proposition}The group $Z_{L}^2(A)$ 
acts by convolution on the right on the set of crossed systems over $R$
corresponding to a fixed measuring $\lu$.  
\end{proposition}
\pf Let $(\lu,\sigma)$ be a crossed system of $A$ on a $k$-algebra
$R$ and let $\omega$ be a lazy $2$-cocycle. Let then
$a,\,b\in A$ and $x\in R$. We have:
$$
\begin{array}{l}
(\a1\lu(\b1\lu x))(\sigma*\omega)(\a2,\b2)\\
\phantom{\sum}=(\a1\lu(\b1\lu x))\sigma(\a2,\b2)\omega(\a3,\b3)\\
\phantom{\sum}=\sigma(\a1,\b1)(\a2\b2\lu x)\omega(\a3,\b3)\\
\phantom{\sum}=\sigma(\a1,\b1)\omega(\a2,\b2)(\a3\b3\lu x)\\
\phantom{\sum}=(\sigma*\omega)(\a1,\b1)(\a2\b2\lu x).
\end{array}
$$
The other relation is proved similarly.
%and
%$$
%\begin{array}{l}
%\sum(\sigma*\omega)(\g1,\h1)(\sigma*\omega)(\g2\h2, m)\\
%\phantom{}=\sum
%\sigma(\g1,\h1)\sigma(\g2\h2,\m1)\omega(\g3,\h3)
%\omega(\g4\h4,\m2)\\
%\phantom{}=\sum
%(\g1\lu\sigma(\h1,\m1))\sigma(\g2,\h2\m2)
%\omega(\h3,\m3)\omega(\g3,\h4\m4)\\
%\phantom{}=\sum
%(\g1\lu\sigma(\h1,\m1)\omega(\h2,\m2))\sigma(\g2,\h3\m3)
%\omega(\g3,\h4\m4)\\
%\phantom{}=\sum
%(\g1\lu((\sigma*\omega)(\h1,\m1))(\sigma*\omega)(\g2,\h2\m2)
%\end{array}
%$$
Hence $(\lu,\,\sigma*\omega)$ is again a crossed system corresponding
to the measuring $\lu$. Then it is clear that
$(\lu,\sigma)\to(\lu,\,\sigma*\omega)$ defines a right action of
$Z^2_{L}(A)$.\hfill$\Box$

\begin{example}{\rm Let $H_4$ be Sweedler's Hopf algebra. 
%Then lazy cocycles
%  are parametrized by $k$ and they are those appearing in \cite{gio},
%  i.e., they are of the form $\sigma_d$ for $d\in k$ such that:
%$$\sigma_d(g, g)=1,$$ 
%$$\sigma_d(g, x)=\sigma_d(x, g)=\sigma_d(g,
%gx)=\sigma_d(gx, g)=0$$
%$$\sigma_d(x, x)=\sigma_d(gx, x)=-\sigma_d(x,
%gx)=-\sigma_d(gx, gx)=\frac{d}{2}.$$.
The action of $Z^2_{L}(A)$ on the set of all crossed systems,
parametrized as in \cite[Table]{[Ma]} by $t$-uples
$(\alpha,\,\delta,\,u,\,a,\,b,\,s)$ is as follows. If $\sigma$
corresponds to $(\alpha,\,\delta,\,u,\,a,\,b,\,s)$ and $\sigma_d$ is
as in Example \ref{swee}, then
$\sigma*\sigma_{d}$ will be the cocycle corresponding to
$(\alpha,\,\delta,\,u,\,a+\frac{d}{2}1,\,b,\,s)$ so the orbits are
parametrized by the $t$-uples $(\alpha,\,\delta,\,u,\,0,\,b,\,s)$.}
\end{example}

After having the notion of a lazy cocycle, we should expect
to have a notion of lazy crossed system. The route for such
a definition is clearly indicated by the straightforward generalization
of lazy Galois objects.

Let $Z$ be a right $A$-comodule algebra. Recall that $R\subset Z$ is
said to be a \textbf{right $A$-Galois extension} if
$Z^{{\rm co} A} = R$ and if 
the linear map
$\kappa_r$ defined by the composition
\begin{equation*}
\begin{CD}
\kappa_r : Z \otimes Z @>1_Z \otimes \rho>>
Z \otimes Z \otimes A @>m_Z \otimes 1_A>> Z \otimes A
\end{CD}
\end{equation*}
induces an isomorphism $Z\otimes_R Z\cong Z\otimes A$.

\begin{defi}
A right $A$-Galois extension $R \subset Z$
is said to be \textbf{lazy} if there exists
a left $R$-linear right $A$-colinear isomorphism
$\psi : R \otimes A \longrightarrow Z$ such that 
$\psi(1\otimes 1) = 1$ and such that the  morphism
$$\beta_\psi := ({\rm id}_A \otimes m_Z) \circ (\tau \otimes {\rm 
id}_Z)
\circ (\psi^{-1} \otimes \psi_{|A}) \circ \rho
: Z \longrightarrow A \otimes Z$$
is an algebra morphism. Such a map $\psi$ is called a \textbf{symmetry 
morphism}
for $R \subset Z$.
\end{defi}

We have the following corresponding definition at the crossed system 
level.

\begin{defi}
A crossed system $(\rightharpoonup,\sigma)$ over $R$
is said to be \textbf{lazy} if
$$a_1b_1 \otimes
((a_2 \rightharpoonup y)\sigma(a_3,b_2)) =
a_3b_2 \otimes ((a_1 \rightharpoonup y) \sigma(a_2,b_1)),
\quad \forall a,b \in A, \ \forall y \in R.$$
\end{defi}

We have the following generalization of Proposition \ref{equivalence}. 
The proof is completely similar and is left to the reader.

\begin{prop}\label{generalized}
Let $R \subset Z$ be a right $A$-Galois extension. Then the
following assertions are equivalent.
\begin{enumerate}
\item $R \subset Z$ is a lazy right $A$-Galois extension.
\item There exists a lazy crossed system
$(\rightharpoonup , \sigma)$ over $R$ such that 
$R \sharp_\sigma \! A \cong Z$ as right $A$-comodule algebras.
\end{enumerate}
\end{prop}

Let $(\rightharpoonup, \sigma)$ be a lazy crossed system over $R$. Then 
the map
$\beta : R \sharp_\sigma A \longrightarrow A \otimes R \sharp_\sigma 
A$,
$x \sharp a \longmapsto
a_1 \otimes x \sharp a_2$, is an algebra morphism
(this is $2  \Rightarrow 1)$ in the proof of Proposition 
\ref{generalized}).
In fact $\beta$ endows $R \sharp_\sigma A$ with a left $A$-comodule
algebra structure and $R \subset Z$ is a left $A$-Galois extension.
This leads to consider general biGalois extension, a notion
which seems not to have been studied before, although 
the definition requires no imagination. 

\begin{defi}
Let $A$ and $B$ be some Hopf algebras. Let $Z$ be an
$A$-$B$-bicomodule algebra. We say that $R \subset Z$
is an $A$-$B$-biGalois extension if
$R \subset Z$ is a left $A$-Galois extension and
is a right $A$-Galois extension.
An $A$-$A$-biGalois extension $R \subset Z$ 
is said to be \textbf{bicleft} if there exists  
a left $R$-linear bicolinear isomorphism
$R \otimes A \cong Z$.
\end{defi}

It is clear that if 
$(\rightharpoonup, \sigma)$ is a lazy crossed system over $R$, then
$R \subset R \sharp_\sigma A$ is a bicleft $A$-$A$-biGalois object.
Similarly to Proposition \ref{bicleft}, we have the following result.

\begin{prop}
Let $R \subset Z$ be an $A$-$A$-biGalois extension. Then the following
assertions are equivalent:
\begin{enumerate}
\item $R\subset Z$ is bicleft.
\item  There exists a lazy crossed system $(\rightharpoonup, \sigma)$
  over $R$ such that 
$R \sharp_\sigma A \cong Z$ as $A$-$A$-bicomodule algebras.
\end{enumerate}
\end{prop}

We conclude with a few words concerning possible generalizations of 
the second lazy cohomology group. 
We would have liked to be able to compose 
lazy crossed systems $(\rightharpoonup, \sigma)$ over $R$ with
$\rightharpoonup$ fixed, in order to have a cohomology
with coefficients in possibly non-commutative algebras.
However we have found that for doing this, one 
needs to require that $\rightharpoonup$ is trivial and $R$ is
commutative. In this case we can define in a straigforward manner
groups $H_L^2(A,R)$, which just turns out to be
$H^2_L(R \otimes A)$ when one works in the category of $R$-algebras.

\end{document}